\documentclass[12pt]{article}
\usepackage{amsxtra,amssymb,amsthm,amsmath,latexsym,verbatim,epsfig}
\usepackage{diagrams}

\theoremstyle{plain}
\newtheorem{definition}{Definition}
\newtheorem{theorem}{Theorem}
\newtheorem{proposition}{Proposition}

\newtheorem{lemma}{Lemma}
\newtheorem{corollary}{Corollary}
\newtheorem{remark}{Remark}
\newtheorem{example}{Example}
\numberwithin{equation}{section}

\newcommand{\refT}[1]{Theorem~\ref{T:#1}}
\newcommand{\refS}[1]{Section~\ref{S:#1}}

\newcommand{\refL}[1]{Lemma~\ref{L:#1}}
\newcommand{\refE}[1]{Example~\ref{L:#1}}
\newcommand{\refC}[1]{Corollary~\ref{C:#1}}

\newcommand{\refP}[1]{Proposition~\ref{P:#1}}
\newcommand{\refD}[1]{Definition~\ref{D:#1}}
\newcommand{\refR}[1]{Remark~\ref{R:#1}}

\def\nd{\noindent}

\def\ve{{\varepsilon}}

\def\lra{\longrightarrow}
\def\hra{\hookrightarrow}

\def\R{{\mathbb R}}
\def\T{{\mathbb T}}

\def\C{{\mathbb C}}
\def\D{{\mathbb D}}
\def\CP{\mathbb{CP}}

\def\calC{{\mathcal C}}

\def\calM{{\mathcal M}}

\def\calO{{\mathcal O}}
\def\hatcalC{{\widehat{\mathcal C}}}

\def\hatg{{\widehat{g}}}
\def\hatw{{\widehat{\omega}}}

\def\dotk{\dot{k}}

\def\dotp{\dot{p}}
\def\dotg{\dot{\gamma}}
\def\Int{\int\limits}
\def\Sum{\mathop\Sigma}

\def\tildeC{{\widetilde{C}}}
\def\tildeJ{{\widetilde{J}}}
\def\tildeL{{\widetilde{L}}}
\def\tildeM{{\widetilde{M}}}
\def\tildeN{{\widetilde{N}}}

\def\tildeX{{\widetilde{X}}}
\def\tildeY{{\widetilde{Y}}}
\def\tildeg{{\widetilde{g}}}
\def\tildeG{{\widetilde{\Gamma}}}
\def\tildeo{{\widetilde{\omega}}}
\def\tildep{{\widetilde{\rho}}}
\def\tildenabla{{\widetilde{\nabla}}}

\def\barD{{\overline{\D}}}
\def\bareta{{\overline{\eta}}}
\def\barxi{{\overline{\xi}}}
\def\bardxi{{\overline{d\xi}}}
\def\bara{{\overline{\alpha}}}
\def\barb{{\overline{\beta}}}

\def\barz{{\overline{z}}}
\def\bardz{{\overline{dz}}}
\def\bareta{{\overline{\eta}}}
\def\bardeta{{\overline{d\eta}}}
\def\bardalpha{{\overline{d\alpha}}}
\def\bardbeta{{\overline{d\beta}}}
\def\bardlambda{{\overline{d\lambda}}}
\def\oH{\buildrel\circ\over H}
\def\oH1{\buildrel\circ\over H\kern-.02in{}^1}
\def\hookuparrow{{\cup\kern-.04in{}^\uparrow}}

\def\l{\ell}

\def\const{\mathop{\rm const}\nolimits}
\def\dist{\mathop{\rm dist}\nolimits}
\def\diam{\mathop{\rm diam}\nolimits}
\def\pr{\mathop{\rm pr}\nolimits}
\def\II{\mathop{\rm II}\nolimits}
\def\CT{{\rm CT}}
\def\Area{{\rm Area}}
\def\id{{\rm id}}
\def\tildeA{{\widetilde{\Area}}}

\def\bee{\begin{equation*}}
\def\eee{\end{equation*}}
\def\be{\begin{equation}}
\def\ee{\end{equation}}

\begin{document}

\title{Conormal bundles to knots and the Gopakumar-Vafa conjecture
\thanks{Math subject classification, Primary: 32Q65, 53D05, 53D12; Secondary: 53D45, 53C21}
\thanks{Keywords: bounded geometry, conifold transition, Gopakumar-Vafa conjecture, Gromov convergence, Gromov-Witten invariants, Lagrangian submanifolds, moduli spaces, pseudoholomorphic curves, symplectic form.}}

\author{Sergiy Koshkin\thanks{partially supported by grant DMS-0204651  from the National Science Foundation.}\\
Department of Mathematics,\\
Northwestern University, Evanston, IL 60208 USA\\
email: koshkin@math.northwestern.edu}

\date{}

\maketitle\thispagestyle{empty}

\begin{abstract}
We offer a new construction of Lagrangian submanifolds for the Gopakumar-Vafa conjecture relating the Chern-Simons theory on the 3-sphere and the Gromov-Witten theory on the resolved conifold. Given a knot in the 3-sphere its conormal bundle is perturbed to disconnect it from the zero section and then pulled through the conifold transition. The construction produces totally real submanifolds of the resolved conifold that are Lagrangian in a perturbed symplectic structure and correspond to knots in a natural and explicit way. We prove that both the resolved conifold and the knot Lagrangians in it have bounded geometry, and that the moduli spaces of holomorphic curves ending on the Lagrangians are compact in the Gromov topology. 
\end{abstract}

\section*{Introduction}\label{S:0}

In \cite{W92} E.\,Witten argues that the large $N$ expansion of the $U(N)$ Chern-Simons theory on a 3-manifold $M$
should be equivalent to an open string theory on $T^\ast M$. In
the absence of knots (Wilson loops) the latter is supposed to describe
pseudoholomorphic curves (strings) on $T^\ast M$ with boundaries on the zero section.
When knots are present H.\,Ooguri and C.Vafa suggested that curves should
additionally be allowed to end on conormal bundles to them.
\cite{OV}. Unfortunately, ordinary Gromov-Witten theory on $T^\ast M$ is
trivial since there are no non-trivial pseudoholomorphic curves there.
Neither closed surfaces nor surfaces ending on the zero section or a
conormal bundle may be pseudoholomorphic due to a vanishing theorem
in \cite{W92} (see \refR{1}). One way around this proposed by E.Witten
himself is to use some degenerate `curves' (fat-graphs) but it is
unclear how to formalize such a theory (see however \cite{Ng}).

Another way around this difficulty was proposed by R.Gopakumar and
C.Vafa in \cite{GV} for $M=S^3$. The idea is to change the topology of
$T^\ast S^3$ so that Gromov-Witten theory on the resulting manifold is
non-trivial and still equivalent to the Chern-Simons theory on $S^3$.
The resulting manifold in this case is the $\calO(-1)\oplus \calO(-1)$
bundle over $\CP^1$ and it can be obtained from $T^\ast S^3$ by
shrinking the $S^3$-cycle to a point and then inserting an
$S^2\hbox{-cycle}$ in its place. Thus, the $U(N)$ Chern-Simons theory on $S^3$ is predicted to be dual to the Gromov-Witten theory on $\calO(-1)\oplus\calO(-1)$. This is the Gopakumar-Vafa conjecture.

The midpoint in the transition is a singular
variety $\calC$ called the conifold and the change from $T^\ast S^3$ to $\calO(-1)\oplus\calO(-1)$ is called the conifold transition. In the physical literature $T^\ast S^3$ is referred to as the deformed
conifold and the resulting $\calO(-1)\oplus\calO(-1)$ bundle $\hatcalC$
as the resolved conifold. Schematically,
$$
\begin{diagram}
T^\ast S^3 & \rDashto  &        &               & \hatcalC \\
           & \rdTo_F   &        & \ldTo_{\pi_2} &          \\ 
           &           &  \calC &               &         
\end{diagram}
$$
with $F$ being the contraction map and $\pi_2$ projecting to $\C^4$ (see \refS{1} for details). We note that $\pi_2^{-1}$ is defined on $\calC\setminus\{0\}$ so the dashed arrow is 'almost' well-defined and $\pi_2^{-1}(0)\simeq\CP^1$ is the exceptional $S^2\hbox{-cycle}$.

As the $S^3$-cycle represented by the zero section in $T^\ast S^3$
shrinks the open curves that end on it become closed
and then get lifted to $\hatcalC$.
Although this picture has no mathematical meaning, note that
unlike $T^\ast S^3$ the resolved conifold $\hatcalC$ does admit non-trivial closed holomorphic curves
and one can talk about equivalence or duality between the Chern-Simons on
$S^3$ and the Gromov-Witten on $\hatcalC$. In the case of closed curves it was verified by a
direct computation in \cite{FP}.

When knots are present the geometric part of the Gopakumar-Vafa conjecture predicts that the conormal bundle
to a knot undergoes the conifold transition and produces some Lagrangian
submanifold $ L$ in $\hatcalC$ \cite{OV}. Then the Chern-Simons theory with knot
observables (Wilson loops) on $S^3$ has to be dual to the open
Gromov-Witten theory where the curves end on $L$.
Ooguri and Vafa were able to produce $L$ explicitly in the case
of the unknot using antiholomorphic involutions. This is a trick
that does not generalize to any other knots. For this case the conjecture
has been verified in \cite{KL,LS} using some narrow definitions
of open Gromov-Witten invariants. 
Later J.\,Labastida, M.\,Mari\~no and C.Vafa offered a way to construct Lagrangians 
for algebraic knots, in particular torus knots \cite{LMV}. This construction as explained by C.H.\,Taubes \cite{T} begins with producing a two-dimensional Lagrangian surface in $\C^2$ that intersects spheres of large radii along the given knot. Then this surface is translated with twisting along the fibers of $\hatcalC$ over the equator of $\CP^1$ completing a half-twist after the full circle (analogous to the M\"obius strip considered as a bundle over the circle). To match up the ends the original surface in $\C^2$ must be centrally symmetric which imposes a restriction on admissible knots. Taubes came up with a generalization to non-algebraic knots and links but the rather artificial symmetry restriction remained. In particular, it excludes something as simple as the trefoil knot. The main flaw of this construction though is that the Lagrangian submanifold constructed is entirely unrelated to the conormal bundle in $T^\ast S^3$ it is supposed to come from. 

Our approach in contrast will be to obtain the corresponding Lagrangians directly  by applying the conifold transition to the conormal bundle $N^\ast_k$ of a knot $k$. As a result the manifold $L$ is produced for all knots in a uniform way and without restrictions. However, this approach presents its own difficulties. Since $N^\ast_k$ intersects the zero-section of $T^\ast S^3$ which is being shrunk into the conifold singularity, it also acquires a singularity in the process. In general this singularity is not resolved by subsequent lifting to $\hatcalC$. The intuitive idea held by physicists (e.g., C.Vafa) is that one needs to perturb $N^\ast_k$ into $N^\ast_{k,\ve}$ disconnected from the zero-section and only then perform the conifold transition to get $\CT(N^\ast_{k,\ve})$. The problem is that $N^\ast_k$ is an exact Lagrangian and there is an obstruction to disconnecting it from the zero section by a symplectic isotopy. It was discovered originally by M.Gromov for compact submanifolds \cite{Gr,ALP} and then generalized to non-compact ones by Y.-G. Oh using the Floer homology \cite{Oh1}.

However, a smooth disconnecting isotopy can easily be found. The construction is very straightforward. Let $k:S^1\hra S^3$ be a naturally parametrised knot. Embed $S^3$ into $\R^4$ in the standard way then the embedding of $T^\ast S^3\simeq TS^3$ into $\R^4\times\R^4$ is also standard. The conormal bundle to $k$ is realized as 
$$
N^\ast_k= \{ (x,p)\in T^\ast S^3\mid x=k(t),\ p\cdot\dotk(t)=0 \}.
$$
To disconnect $N^\ast_k$ from the zero-section we shift it in the fibers of $T^\ast S^3$ in the direction tangent to the knot, namely
$$
N^\ast_{k,\ve}:=\{(x,p+\ve\dotk(t))\mid x=k(t),\ p\cdot\dotk(t)=0\}. 
$$  
Since $F(N^\ast_{k,\ve})$ misses the conifold singularity the conifold transition is given simply by 
$$
\CT(N^\ast_{k,\ve}):=\pi_2^{-1}\circ F(N^\ast_{k,\ve}).
$$
The transition $\CT(N^\ast_{k,\ve})$ turns out to be a smooth submanifold in $\hatcalC$ and has the correct topology $S^1\times\R^2$. Predictably, it fails to be Lagrangian in the standard K\"ahler structure of the resolved conifold. 

We beleive however that this is a false problem. Although the standard K\"ahler structure on $\hatcalC$ is the simplest one it is not in any way special from the physical point of view. The physically significant structure, if any, is the one induced by the Calabi-Yau metric. The Calabi-Yau metric on $\hatcalC$ is known almost explicitly \cite{CdO} but no attempt has been made to check if even the Ooguri-Vafa submanifold for the unknot is Lagrangian in it. Moreover, the computations of open invariants in \cite{KL} only use the fact that it is Lagrangian in the standard metric on the resolved conifold. On the other hand, one does not need a necessarily Lagrangian submanifold to build a theory of open holomorphic curves. It suffices to have a totally real submanifold \cite{Oh} with some uniformity conditions in non-compact cases \cite{S,ALP}. Conifold transitions of perturbed conormal bundles constructed in this paper do meet these conditions. Moreover, one can show that the moduli of holomorphic curves ending on these submanifolds are compact and thus suitable for defining open Gromov-Witten invariants (\refT{5}).  

The paper is organized as follows. In \refS{1} we briefly review the conifold transition and introduce a natural notion of the conifold transition for submanifolds of $T^\ast S^3$. In \refS{2} we compute explicitly the conifold transitions of the conormal bundles to the unknot and torus knots. In the former case we get the well-known Ooguri-Vafa Lagrangian \cite{OV}, while in the latter a variety which is neither smooth nor Lagrangian. In \refS{3} a perturbed conormal bundle is defined and we prove that its conifold transition is a tame Lagrangian in $\hatcalC$ (see \refD{5}). In \refS{4} we lay the geometric groundwork for the compactness result in \refS{5}. Namely, we use the technique of second fundamental forms and bi-Lipschitz maps to prove that the resolved conifold has bounded geometry, i.e. its sectional curvature is bounded from above and the injectivity radius is bounded from below. The key point is \refL{10} which formalizes the idea that the conifold has cone-like geometry. Finally, in \refS{5} the moduli spaces of open curves are introduced following \cite{L} and the compactness of the moduli of curves ending on $\CT(N^\ast_{k,\ve})$ is proved. In the end we present our conclusions.

\smallskip

{\bf Acknowledgements.} The author would like to thank D.Auckly without whose advice and encouragement this paper would never have been completed. I am also grateful to M.Liu for patiently explaining to me the theory of open Gromov-Witten invariants and C.Vafa and M.Mari\~no for suggesting the construction of perturbed conormal bundles. Main ideas of this paper were discussed during the program Symplectic Geometry and Physics at the Institute  of Pure and Applied Mathematics in March-June 2003 and the workshop Interaction of finite type and Gromov-Witten invariants at the Banff International Research Station in November 2003. It is my pleasant duty to thank all their organizers and participants.

\section{The conifold transition.}\label{S:1}

In this section we first review some basic facts about the conifold transition \cite{Cl,CdO,GR,ST} and fix the notation used throughout the paper. Then we define the conifold transition for submanifolds in $T^\ast S^3$ and apply it to conormal bundles to knots in $S^3$.

It is convenient to think of $S^3$ as being embedded into $\R^4$ as the unit sphere and identify $T^\ast S^3$ with $TS^3$ via the standard metric. At each point $x$ of the sphere the tangent space $T_xS^3$ is naturally identified with the tangent hyperplane in $\R^4$ at this point. Shifting it to the origin we get the  subspace of $\R^4$ orthogonal to $x$ and obtain a natural realization of the tangent bundle in $\R^4\times\R^4$ as 
\bee
T^\ast S^3=TS^3=\{(x,p)\in \R^4\times \R^4\mid |x|=1,\,p\cdot x=0\}.
\eee 
Now introduce complex coordinates on $\R^4\times \R^4\simeq\C^4$ by $z_j:=x_j+ip_j$. As realized above $T^\ast S^3$ is not an algebraic submanifold of $\C^4$ but it is diffeomorphic to any member of the family with $a>0$:
\bee \calC_a:=\{z \in \C^4\mid \Sum_j z^2_j=a^2\} 
   = \{(x,p)\in \R^4\times \R^4\mid  |x|^2-|p|^2=a^2, p\cdot x =0\}.
\eee
Any $\calC_a$ is manifestly algebraic and the diffeomorphism with $T^\ast S^3$ is given by
\bee
  \begin{aligned}
  T^\ast S^3   & \overset{F_a}{\lra} \calC_a\\
  (x,p)        &\longmapsto (x\sqrt{a^2+|p|^2},p).
 \end{aligned}
\eee
The standard K\"ahler form 
\be\label{e1.001}
\omega:=\frac{i}{2} \sum_j dz_j\wedge\overline{dz_j}=\sum_jdx_j\wedge dp_j
\ee
on $\C^4$ restricts to a K\"ahler form on every $\calC_a$. We can pull these forms back to $T^\ast S^3$ via $F_a$, i.e. $\omega_a:=F^\ast_a\omega\left|_{\calC_a}\right.$. It is straightforward to check that the standard symplectic form on $T^\ast S^3$ is obtained as the limit $\omega_\infty:=\lim_{a\to\infty}\frac{1}{a}\omega_a$ but $\omega_\infty$ is no longer K\"ahler.

\begin{figure}\label{F:1}
\hskip1in\epsfig{file=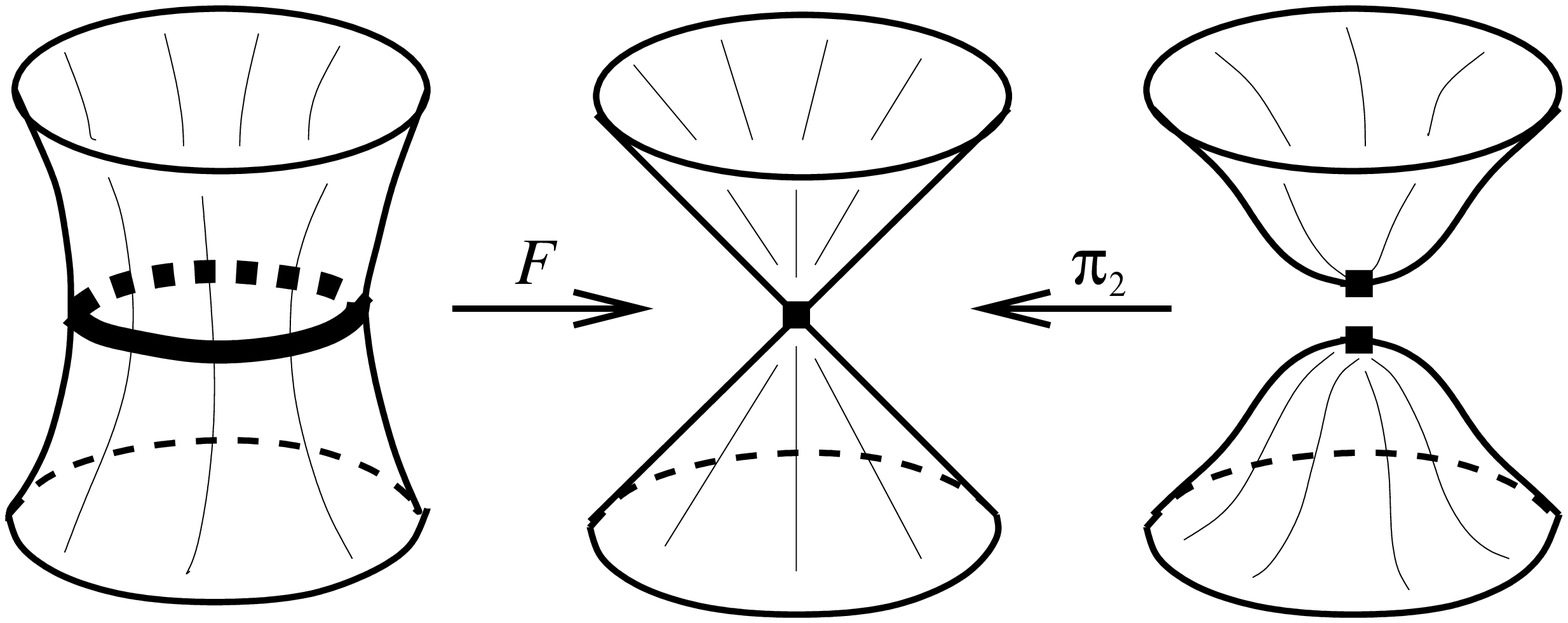,width=4truein}
\caption{The 'conifold' transition two dimensions down $S^1\times\R^1\leadsto S^0\times\R^2$.}
\end{figure}

Let us look at the other limit $a\to0$ (cf. \cite{ST}, sec.3). The algebraic subvariety
\be\label{e1.1}
\calC:=\calC_0=\{z \in\C^4\mid \Sum_jz^2_j=0\}
 = \{(x,p)\in \R^4\times\R^4|\ |x|=|p|,p\cdot x=0\}
\ee
is called the {\bf conifold} \cite{CdO}. It is singular with a nodal (ordinary double) point at the origin \cite{Cl,GR}. 
Equivalently, one can think of it as $T^\ast S^3$ with the singular form $\omega_0:=F^\ast_0\omega\left|_{\calC}\right.$ that degenerates along the zero section. The manifolds $\calC_a$ or equivalently
$(T^\ast S^3, \omega_a)$ are called the {\bf deformed conifolds}. The parameter $a$ has a simple geometric interpretation as the radius of the 3-sphere $p=0$ in $\calC_a$, i.e., of the zero section. As $a$ goes to $0$ this
sphere collapses into the singular point in $\calC$. The map 
\be\label{e1.02}
F(x,p):=F_0(x,p)=(|p|x,p)
\ee contracts $T^\ast S^3$ onto the conifold (Fig.\,1).

An alternative to smooth deformation when dealing with singularities is resolution. The simplest type of resolution is blow-up \cite{MS} and blowing up the origin in $\C^4$ produces an algebraic submanifold of $\CP^3\times \C^4$. Let $\pi_2:\CP^3\times\C^4\to\C^4$ be the natural projection to the second factor. Then the proper transform 
$\widetilde{\calC}:=\overline{\pi^{-1}_2 (\calC\setminus\{0\})}$ is the blow-up or the large resolution of $\calC$ (here and below the overline denotes the closure). The inverse image of the singular point $\pi^{-1}_2(0)\simeq \CP^2$ is called the exceptional divisor. 

However, the conifold singularity admits a smaller resolution that only adds an exceptional curve $\CP^1$ instead of a whole divisor $\CP^2$. To describe this resolution it is convenient to use different coordinates on $\C^4$:
\be\label{e1.01}
\begin{matrix} w_1 &=z_1+iz_2 & \qquad & w_2&=-z_3+iz_4\\
                  w_3 &=z_3+iz_4 & \qquad & w_4&=z_1-iz_2. 
\end{matrix}
\ee
Up to the factor of $\frac{1}{\sqrt{2}}$ this is a unitary transformation of
$\C^4$ and the difining equation (\ref{e1.1}) of the conifold becomes:
\bee
\calC=\left\{w\in\C^4\mid w_1w_4-w_2w_3=\begin{vmatrix}w_1 & w_2\\w_3&w_4\end{vmatrix}=0\right\}.
\eee
In the $w$-coordinates the small resolution can be written explicitly as
\be\label{e1.2}
\hatcalC :=\left\{
   ([u:v], w_1,w_2,w_3,w_4)\in \CP^1\times\C^4
   \left|
     \begin{vmatrix}w_1&w_2\\u&v\end{vmatrix}
    =\begin{vmatrix}w_3&w_4\\u&v\end{vmatrix}
   \right. =0
  \right\},
\ee
and $\hatcalC$ is called the {\bf resolved conifold}. It is obvious from the definition that 
$\hatcalC\to\CP^1$ is a holomorphic subbundle of the trivial
$\C^4\hbox{-bundle}$ over $\CP^1$, and therefore the total space is a
smooth manifold. We denote its {\bf zero section} by $0(\hatcalC)$. The resolution preserves the canonical class \cite{Cl,GR} and $\hatcalC$ along with $\calC$ and $\calC_a$ is a Calabi-Yau threefold \cite{J}.

To understand the resolved conifold better consider the tautological line bundle over $\CP^1$:
$$
\calO(-1):=\left\{([\lambda],w)\in\CP^1\times\C^2\mid w\in[\lambda]\right\}.
$$
Here we interpret $[\lambda]$ as a line in $\C^2$ and the fiber over it consists of all the points on this line (hence the name tautological). The letter $\calO$ is traditional for holomorphic bundles in algebraic geometry and and the number $-1$ refers to the fact that the first Chern class of this bundle evaluates to $-1$ on the base $\CP^1$ \cite{BT,MS}. More explicitly,
\bee
 \begin{aligned}
   \calO(-1)=\left\{
   ([u:v], w_1,w_2)\in \CP^1\times\C^2
   \left|
     \begin{vmatrix}w_1&w_2\\u&v\end{vmatrix}
   \right. =0
  \right\}
\end{aligned}
\eee
and one can see by inspection that
\bee
\hatcalC=\calO(-1)\oplus\calO(-1).
\eee
Denoting by $\pi_1$, $\pi_2$ the natural projections to $\CP^1$ and
$\C^4$ respectively, one notes that 
$$
\calC=\pi_2(\hatcalC),\quad \hatcalC=\overline{\pi^{-1}_2(\calC\setminus\{0\})}.
$$ 
Moreover, the projection $\pi_2$ restricts to a biholomorphism from $\hatcalC\setminus \pi^{-1}_2(0)$ to $\calC\setminus\{0\}$. The resolved conifold $\hatcalC$ admits the K\"ahler form $\pi^\ast_1 \omega_{FS}+\pi^\ast_2\omega$, where $\omega_{FS}$ is the Fubini-Study form on $\CP^1$ and $\omega$ is the standard K\"ahler form \eqref{e1.001} on $\C^4$.

The conifold transition from the deformed to the resolved conifold is
streamlined in the following diagram
\be\label{e1.4} 
\begin{diagram}
T^\ast S^3 & \rDashto  &        &               & \hatcalC & \rInto &\CP^1\times\C^4\\
           & \rdTo_F   &        & \ldTo_{\pi_2} &          &        &        \\ 
           &           &  \calC &               &          &        &
\end{diagram}
\ee 
where $F(x,p)=(|p|x,p)$ is the contraction \eqref{e1.02} that shrinks the zero section to
the conifold singularity at the origin. Since $\pi^{-1}_1(0)\simeq\CP^1$
the singularity gets replaced by a new $S^2\hbox{-cycle}$ in $\hatcalC$. Topologically we have the transition $S^3\times\R^3\leadsto S^2\times\R^4$ (see Fig.\,1).

We will also be interested in the conifold transitions of certain submanifolds of $T^\ast S^3$. The most natural transition seems to be the contraction by $F$ into $\calC$ followed by the proper transform. 
\begin{definition}\label{D:0}
The {\bf conifold transition of a submanifold} $N\subset T^\ast S^3$ is
\be\label{e1.04} 
\CT(N):=\overline{\pi^{-1}_2(F(N)\setminus\{0\})}.
\ee 
\end{definition}
In our case $F(N)$ will be Lagrangian and not complex
submanifolds so we can not expect that the proper transform will
produce smooth manifolds. In fact, as examples in the next section show,
the conifold transitions are not smooth in general.

The submanifolds we are primarily interested in are conormal bundles to knots. We will show in \refE{1} that our conifold transition of the conormal bundle to the unknot is the Ooguri-Vafa Lagrangian obtained in \cite{KL,OV} as the fixed locus of an antiholomorphic involution.
\begin{definition}\label{D:1}
Let $S\hra M$ be a submanifold. Then its {\bf conormal bundle} in
$T^\ast M$ is
\bee
 N^\ast_S:=\{\l \in T^\ast M\mid \pi(\l)\in S,\ \l\mid_{TS}=0\}
\eee
\end{definition}
Let $k:S^1\to S^3$ be a knot. Then under the identification
of $T^\ast S^3$ with the submanifold of $\R^4\times\R^4$ one gets:
\be\label{e1.5}
  \begin{aligned}
  N^\ast_k
   & =\{(k(t),p)\in T^\ast S^3 \mid t\in S^1,\ p\cdot\dotk(t)=0\} \\
   & =\{(k(t),p)\in \R^4\times\R^4\mid t\in S^1,\ p\cdot k(t)=p\cdot\dotk(t)=0\}.
   \end{aligned}
\ee
\begin{lemma}\label{L:0}
$N^\ast_k$ admits the following parametrization
\be\label{e1.6}
  \begin{aligned}
  S^1\times \R^2   & \overset{K}{\lra} T^\ast S^3\\
  (t,\alpha,\beta) &\longmapsto (k(t),\alpha p^1(t)+\beta p^2(t))
 \end{aligned}
\ee
where $p^1(t)$, $p^2(t)$ are fundamental solutions in $\R^4$ to
\bee
  \begin{cases} k(t)\cdot p &=0\\ \dotk(t)\cdot p& =0,
  \end{cases}
\eee
and for every $t$ the vectors $k(t)$, $\dotk(t)$, $p^1(t)$, $p^2(t)$
form an orthonormal basis in
$\R^4$. Moreover, $p^1$, $p^2$ are $C^\infty$-smooth if $k$ is.
\end{lemma}

\begin{proof}
Since $k(S^1)\subset S^3$ we have $|k(t)|=1$ which implies
$k(t)\cdot\dotk(t)=0$.
Choosing the natural parametrization for the knot we also get
$|\dotk(t)|=1$. This means that the system for $p^1$, $p^2$ is
non-degenerate. Choose a stereographic projection
$\sigma:S^3\to\R^3\cup\{\infty\}$ so that $k$ passes neither through the
north nor through the south pole. By transversality we may also
assume that $\sigma(k(t))$ and $\sigma_\ast\dotk(t)$ are linearly independent
for every $t$. Set $\tilde{p}^1(t)$ to be the cross-product $\sigma(k(t))\times \sigma_\ast\dotk(t)$ in $\R^3$
then $k(t)$, $\dotk(t)$, $\sigma^{-1}_\ast \tilde{p}^1(t)$ are
linearly independent in $\R^4$. Now one can get $p^1(t)$  by the
Gram-Schmidt process. Finally, $p^2(t):=k(t)\times \dotk(t)\times p^1(t)$,
the cross-product of three vectors in $\R^4$.
The smoothness is obvious from the construction.
\end{proof}

It turns out that $N^\ast_k$ is an exact Lagrangian submanifold in
$T^\ast S^3$. Let us recall the definition \cite{ALP}:

\begin{definition}\label{D:2}
A symplectic manifold $(X,\omega)$ is exact
if $\omega$ has a primitive, i.e. there is a 1-form $\lambda$ such that $\omega=d\lambda$.
A Lagrangian submanifold $L\hra X$ is exact if
$[\lambda\!\mid_L]=0\in H^1(L,\R)$.
\end{definition}

In particular, any $X=T^\ast M$ with the canonical symplectic structure is
exact and $\lambda$ is the canonical Liouville form. In Darboux coordinates
$\lambda=-\Sigma_j p_j dx_j$ and
$\lambda\mid_{0(T^\ast M)}=0$, where $0(T^\ast M)$
denotes the zero section corresponding to $p=0$.
So the zero section is always an exact Lagrangian submanifold.

\begin{remark}\label{R:1}
Note that if $f:(\Sigma,\partial\Sigma)\to (T^\ast M,L)$ is a pseudoholomorphic open curve (see \refS{5}) ending on the zero section or any other exact Lagrangian submanifold $L$ then
\bee\label{e1.13} 
\Area(f)=\int_\Sigma f^\ast\omega=\int_\Sigma d(f^\ast \lambda) =\int_{\partial\Sigma} f^\ast\lambda=[\lambda\!\mid_L](f_\ast[\partial\Sigma])=0. 
\eee
Thus $f$ has to be a constant. E.Witten \cite{W92} calls this fact 'the vanishing theorem'.
\end{remark}

In our case the symplectic form on $T^\ast S^3$ was the restriction of
$\omega=\Sigma^4_{j=1} dx_j\wedge dp_j$ from $\R^4\times \R^4$.
If one sets $\lambda:=-\Sigma^4_{j=1} p_j dx_j$, then obviously
$\omega=d\lambda$ and since $d$ commutes with restrictions
$\omega\mid_{T^\ast S^3}=d(\lambda\!\mid_{T^\ast S^3})$.

\begin{lemma}\label{L:1}
$N^\ast_k$ is an exact Lagrangian submanifold in $T^\ast S^3$.
\end{lemma}

\begin{proof}
In fact, we will show that $\lambda\!\mid_{N^\ast_k}=0$.
Identifying $T_{(x,p)}N^\ast_k$ with a subspace of $\R^4\times\R^4$ in
the usual way one gets:
 \bee\begin{aligned}
  K_\ast\partial_t & =(\dotk,\alpha\dot{p}^1+\beta\dot{p}^2)\\
  K_\ast\partial_\alpha & =(0,p^1)\\
  K_\ast\partial_\beta  & =(0,p^2)\\
  \end{aligned}\eee
Thus
\bee\label{e1.14}
  \lambda(K_\ast\partial_t)=(\alpha p^1+\beta p^2)\cdot\dotk
      =\alpha p^1\cdot\dotk +\beta p^2\cdot\dotk=0\eee
and
\bee\label{e1.15}
  \lambda(K_\ast\partial_\alpha)=(\alpha p^1+\beta p^2)\cdot 0
      = 0 =\lambda(K_\ast \partial_\beta).\eee
Finally, $N^\ast_k$ has the right dimension:
$\dim N^\ast_k=\frac{1}{2}\dim T^\ast S^3=3$.
\end{proof}

Intuitively, this corresponds to the fact that the only
$x\hbox{-direction}$ in $N^\ast_k$ is the one orthogonal
to its $p\hbox{-directions}$ so $\lambda=-p\cdot dx$ vanishes on it
`by definition'.

Recall from \eqref{e1.02} that the first half of the conifold transition is the
contraction $F(x,p)=(|p|x,p)$. It turns out that although this map
is not symplectic it does map exact Lagrangians into exact Lagrangians
away from the singular locus.

\begin{lemma}\label{L:2} 
$F(N^\ast_k)\setminus\{0\}$ is a Lagrangian submanifold in $\calC$.
\end{lemma}

\begin{proof}
On $T^\ast S^3$ one has $F^\ast\lambda=|p|\lambda$, $F^\ast \omega=|p|\omega+d|p|\wedge\lambda$. Indeed,
\bee
  \begin{aligned}
   F^\ast \lambda &=F^\ast(-\Sum_j p_j dx_j)= -\Sum_j p_j d(|p|x_j)\\
   &=-|p|\Sum_j p_j dx_j - \Sum_j p_j x_j d|p|\\
   &=|p|\lambda -(p\cdot x)d|p|=|p|\lambda
   \end{aligned}\eee
since $p\cdot x=0$ on $T^\ast S^3$. The second relation follows from
the first one:
\bee
\begin{aligned}
   F^\ast \omega &=F^\ast d\lambda=d F^\ast \lambda=d(|p|\lambda)\\
            &=|p|d\lambda +d|p|\wedge\lambda
   \end{aligned} 
\eee
$N^\ast_k\cap 0(T^\ast S^3)$ will be mapped into the conifold
singularity at the origin but away from that point
$\lambda\!\mid_{F(N^\ast_k)}=F^\ast\lambda\!\mid_{N^\ast_k}
  =|p|\lambda\mid_{N^\ast_k}=0$,
i.e., $F(N^\ast_k)\setminus\{0\}$ is still a Lagrangian submanifold.
\end{proof}
The last step is to lift $F(N^\ast_k)$ to the resolved conifold
$\hatcalC$. The K\"ahler structure on $\hatcalC$ is induced by the product structure
on $\CP^1\times \C^4$, namely
$\hatw:=(\pi^\ast_1 \omega_{FS}+\pi^\ast_2\omega)\!\mid_\hatcalC$ with
$\omega=\frac{i}{2}\Sigma^4_{j=1} dz_j\wedge\overline{dz_j}
  =i\Sigma^4_{j=1}dw_j\wedge\overline{dw_j}$
in complex coordinates.

\begin{theorem}\label{T:1}
$\pi^{-1}_2(F(N^\ast_k)\setminus\{0\})$ is Lagrangian in $(\hatcalC,\hatw)$
if and only if it projects to a set of zero volume in $\CP^1$.
\end{theorem}

\begin{proof}
Since $\pi^\ast_2:\hatcalC\setminus\pi^{-1}_2(0)\to\calC\setminus\{0\}$
is a biholomorphism we have
\bee\label{e1.18}
\begin{aligned}
 \hatw\mid_{\pi^{-1}_2(F(N^\ast_k)\setminus\{0\})}
 & =\pi^{-1\ast}_2\hatw\mid_{F(N^\ast_k)\setminus\{0\}}
   =\pi^{-1\ast}_2
 (\pi^\ast_1\omega_{FS}+\pi^\ast_2\omega)\mid_{F(N^\ast_k)\setminus\{0\}}\\
 &=(\pi_1\circ \pi^{-1}_2)^\ast\omega_{FS}\mid_{F(N^\ast_k)\setminus\{0\}}
     +\omega\mid_{F(N^\ast_k)\setminus\{0\}}\\
  &=\omega_{FS}\mid_{\pi_1(\pi^{-1}_2(F(N^\ast_k)\setminus\{0\}))},
 \end{aligned}\eee
where we used the fact that $F(N^\ast_k)\setminus\{0\}$ is Lagrangian by \refL{2}.
The left-hand side is $0$ if and only if $\pi^{-1}_2(F(N^\ast_k)\setminus\{0\})$
is Lagrangian, while the right-hand side is $0$ if and only if its projection
to $\CP^1$ has $0$ volume.
\end{proof}

Recall that by our definition
$\CT(N^\ast_k)=\overline{\pi^{-1}_2(F(N^\ast_k)\setminus\{0\})}$
so this theorem does not guarantee that the conifold transition of a
conormal bundle is a manifold even if it does project to a null set.
However, if the closure is indeed a smooth manifold with the projection of zero volume,
it will automatically be Lagrangian ($\hatw=0$ on the closure by continuity).
Another remark is that our choice of $\hatw$ on $\hatcalC$ is more or less
arbitrary. From the physical point of view a more natural choice is
$\omega_{CY}$, the K\"ahler form induced by the Calabi-Yau metric on
$\hatcalC$ \cite{CdO}. But as we will see, $\CT(N^\ast_k)$ is not
even a smooth manifold already for torus knots.

\section{The unknot and torus knots}\label{S:2}

Here we use the parametrization of $N^\ast_k$ from the previous
section (\refL{0}) to compute $\CT(N^\ast_k)$ for $k$  the
unknot or a torus knot. For this computation it is convenient to use
a different complex structure on $\R^4\times\R^4$ given by the new holomorphic coordinates
\bee
  \begin{aligned}
  \xi & =(\xi_1,\xi_2) =(x_1+ix_2, x_3+ix_4)\\
  \eta& =(\eta_1,\eta_2) =(p_1+ip_2,p_3+ip_4).
  \end{aligned}
  \eee
In these coordinates
\bee
  T^\ast S^3=
  \{(\xi,\eta)\in \C^2\times\C^2\mid |\xi|=1,\  Re(\xi\bareta)=0\} \eee
and the change to $w\hbox{-coordinates}$ is
\bee
  \begin{matrix}
  w_1=\xi_1+i\eta_1 & \qquad w_2=-(\barxi_2+i\bareta_2)\\
  w_3=\xi_2+i\eta_2  & \qquad w_4=\barxi_1+i\bareta_1
  \end{matrix} \eee
The parametrization of a conormal bundle has the same form as before
\bee
  K(t,\alpha,\beta)=(k(t),\alpha p^1(t)+\beta p^2(t)) 
\eee
but $k$, $p^1$, $p^2$ now are $\C^2$ vectors.
Applying the contraction $F$ one obtains
\bee\label{e2.5}
  F\circ K(t,\alpha,\beta)
  =(k(t)\sqrt{\alpha^2+\beta^2}, \quad \alpha p^1(t)+\beta p^2(t)), \eee
which suggests the change to polar coordinates
\bee\label{e2.6}
  r=\sqrt{\alpha^2+\beta^2}, \quad \tan\theta=\frac{\beta}{\alpha} \eee
where $F\circ K(t,r,\theta)=r(k(t),\wp(t,\theta))$,
$\wp(t,\theta):=p^1(t)\cos\theta+p^2(t)\sin\theta$.
Here $r\geq 0$ and $r=0$ corresponds to the conifold singularity.

\begin{figure}\label{F:2}
\hskip2.3in\epsfig{file=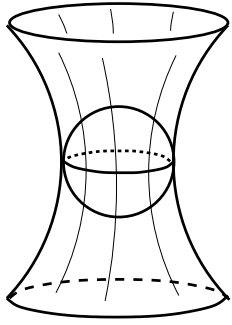,width=1.2truein}
\caption{Conifold transition for the unknot (low-dimensional analog).}
\end{figure}

\begin{example}\label{E:1}
For the unknot we have
\bee\label{e2.7}
  \begin{aligned}
      k(t) &= (e^{it},0)\in S^3\hra \C^2\\
  \dotk(t) &= (ie^{it},0) \\
    p^1(t) &= (0,1) \\
    p^2(t) &= (0,i).
  \end{aligned}
 \eee
 \end{example}
Thus $\wp(t,\theta)=(0,1)\cos\theta+(0,i)\sin\theta =(0,e^{i\theta})$
and $F\circ K(t,r,\theta)=r(e^{it},e^{i\theta})$.
In terms of this parametrization
$F(N^\ast_k)\setminus\{0\}=F\circ K(S^1\times\R_{>0}\times S^1)$.
To find the proper transform we change to $w\hbox{-coordinates}$:
\bee
  \begin{aligned}
  w_1=re^{it} &\qquad w_2=\overline{ire^{i\theta}}=-ire^{-i\theta} \\
  w_3=ire^{i\theta} & \qquad w_4=\overline{re^{it}}=re^{-it}
  \end{aligned}\eee
By definition of the resolved conifold
\bee
  \pi^{-1}_2(w_1,w_2,w_3,w_4)=([u:v],w_1,w_2,w_3,w_4),\eee
where
$\begin{vmatrix}w_1&w_2\\u&v\end{vmatrix}
 =\begin{vmatrix}w_3&w_4\\u&v\end{vmatrix}=0.$
Since on $F(N^\ast_k)\setminus\{0\}$ when $r>0$
both $w_1$, $w_2$ are never $0$
we can just set $[u:v]=[w_1:w_2]$ and get:
\bee
\begin{aligned}
 \CT(N^\ast_k)
 &=\overline{\pi^{-1}_2(F(N^\ast_k)\setminus\{0\})}\\
 &=\{([e^{it}:-ie^{-i\theta}],re^{it},-ire^{-i\theta},ire^{i\theta},re^{-it})
   \mid t,\theta\in S^1,\ r\geq 0\}\\
 &=\{([ie^{i(t+\theta)}:1],re^{it},-ire^{i\theta},ire^{i\theta},re^{-it})
   \mid t,\theta\in S^1,r\geq 0\}
\end{aligned} \eee
This becomes more transparent if one sets
\bee\label{e2.11}
  \begin{aligned}
  \alpha:=&ie^{i(t+\theta)}\\
  b:=&-ire^{-i\theta}
  \end{aligned}\eee
so that
\bee
\CT(N^\ast_k)=\{([\alpha:1],\alpha b,b,\overline{b},\overline{\alpha b})
\mid\alpha\in S^1,b\in\C\}.
\eee  

This is a smooth submanifold of $\hatcalC$ diffeomorphic to
$S^1\times\C$ and since $|\alpha|=1$ it fibers over the equator of
$\CP^1$ (Fig.\,2). By \refT{1} this means that $\CT(N^\ast_k)$ is also
Lagrangian. In fact, this is the same Lagrangian submanifold that was obtained
in \cite{OV} as the fixed locus of an antiholomorphic involution and used in
\cite{KL} to compute open Gromov-Witten invariants. Note that the
topologies of $N^\ast_k$ and $\CT(N^\ast_k)$ are the same,
namely $S^1\times\R^2$ even though $T^\ast S^3$ changes its topology from $S^3\times\R^3$ to $S^2\times\R^4$.

\begin{example}[Torus knots]\label{E:2}
\end{example}
There is a standard copy of a 2-torus sitting in $\C^2:$
$\{(\xi_1,\xi_2)\mid |\xi_1|=|\xi_2|=1\}$.
If we change the normalization from $1$ to $\frac{1}{\sqrt{2}}$
this torus will sit inside $S^3\hra\C^2$.
The embedding $k(t)=\frac{1}{\sqrt{2}}(e^{imt},e^{int})$
obviously winds $m$ times around one of the cycles in $\T^2$ and
$n$ times around the other one. Therefore, for relatively prime $(m,n)=1$ it
represents an $(m,n)$ torus knot.
We will assume $m\not= n$ since $m=n=1$ is the case of the unknot and
otherwise $m,n$ can not be relatively prime. We have
\bee
  \begin{aligned}
  k(t) &= \frac{1}{\sqrt{2}} (e^{imt},e^{int}) \\
  \dotk(t) &= \frac{1}{\sqrt{2}} (ime^{imt},ine^{int}) \\
  p^1(t) &= \frac{1}{\sqrt{2}} (e^{imt}, -e^{int})\\
  p^2(t) &= \frac{1}{\sqrt{m^2+n^2}} (ine^{imt},-ime^{int})
  \end{aligned} \eee
The parameter $t$ here is obviously not the arclength but the only effect this
has is that $|\dotk|^2=\frac{m^2+n^2}{2}=\const$ instead of $1$ so
the difference is insignificant.
\bee
 \wp(t,\theta)=\frac{1}{\sqrt{2}} \left( e^{imt}
 \left(\cos\theta+i\frac{n\sqrt{2}}{\sqrt{n^2+m^2}}\sin\theta\right), -e^{int}
 \left(\cos\theta+i\frac{m\sqrt{2}}{\sqrt{n^2+m^2}}\sin\theta\right)
 \right) \eee
Writing $F\circ K(t,r,\theta)=r(k(t),\wp(t,\theta))$ and
changing to $w\hbox{-coordinates}$ one finds
\bee
  \begin{aligned}
  w_1&=\frac{re^{imt}}{\sqrt{2}}
   \left(1-\frac{n\sqrt{2}}{\sqrt{n^2+m^2}}\sin\theta+i\cos\theta\right)\\
  w_2&=-\frac{re^{-int}}{\sqrt{2}}
   \left(1+\frac{m\sqrt{2}}{\sqrt{n^2+m^2}}\sin\theta-i\cos\theta\right)\\
  w_3&=\frac{re^{int}}{\sqrt{2}}
   \left(1-\frac{m\sqrt{2}}{\sqrt{n^2+m^2}}\sin\theta-i\cos\theta\right)\\
  w_4&=\frac{re^{-imt}}{\sqrt{2}}
   \left(1+\frac{n\sqrt{2}}{\sqrt{n^2+m^2}}\sin\theta+i\cos\theta\right)\\
  \end{aligned}\eee
Just as in the case of the unknot for $r>0$ one can set $[u:v]=[w_1:w_2]=$
\bee
=\left[e^{imt} \left(1-\frac{n\sqrt{2}}{\sqrt{n^2+m^2}}
\sin\theta+i\cos\theta\right) :-e^{-int}  \left(1+\frac{m\sqrt{2}}{\sqrt{n^2+m^2}}
\sin\theta-i\cos\theta\right)\right]
\eee
and since the last expression does not depend on $r$ taking the closure is
simply allowing $r=0$ in the formulas for the $w_j$ above.
Since $F$ and $\pi^{-1}_2$ are diffeomorphisms away from the zero section
of $T^\ast S^3$ and the origin respectively, $\CT(N^\ast_k)$ is a smooth
manifold everywhere but possibly at the zero section of
$\calO(-1)\oplus\calO(-1)\simeq\hatcalC$ where $r=0$.

\begin{figure}\label{F:3}
\hskip2.3in\epsfig{file=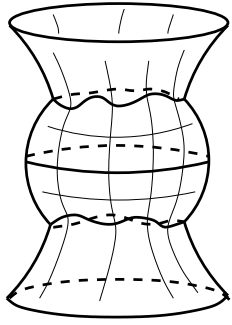,width=1.2truein}
\caption{Corner singularity in the conifold transition for torus knots.}
\end{figure}

Now we take a look at the $\CP^1$ projection of $\CT(N^\ast_k)$.
It is parametrized by $[u:v](t,\theta)$ given above and
$\frac{u}{v}$ is a possible coordinate on the projection (unless $m=n$,                  
which is the case we excluded). Hence
\bee
  w(t,\theta):=\frac{u}{v}=-e^{i(m+n)t}
  \ \frac{1-\frac{n\sqrt{2}}{\sqrt{n^2+m^2}}\sin\theta+i\cos\theta}
       {1+\frac{m\sqrt{2}}{\sqrt{n^2+m^2}}\sin\theta-i\cos\theta}
\eee
parametrizes the image of the projection in $\C$. The first factor is $1$ in absolute value while the absolute value of the second one changes between two positive values less than and greater than $1$. 

This means that the $\CP^1$ trace of $\CT(N^\ast_k)$ is an annulus containing the equator. In the case of the unknot the trace was just the equator itself, in particular it was one-dimensional (Fig.\,2). This is due to the fact that the unknot in $S^3$ can be flat, i.e. placed within a 2-plane, which is impossible for any non-trivial knot. As a result of replacing the 3-cycle by a 2-cycle the conormal bundle to a non-planar knot 'smashes' into an annulus on $\CP^1$ and $\CT(N^\ast_k)$ acquires a corner singularity along the edge of the annulus (see Fig.\,3). By \refT{1} $\CT(N^\ast_k)$ is not Lagrangian even away from the singularity. For the unknot the flatness makes it possible for the conifold transited conormal bundle to just touch the equator.

\section{Perturbed conormal bundles}\label{S:3}

The obvious reason for the conifold transition of a conormal bundle
to be singular is that it intersects the zero-section of $T^\ast S^3$
which collapses into a singular point. The simplest way to avoid this
is to perturb $N^\ast_k$ so that it is disconnected from the zero
section. Of course, we would like to obtain an exact Lagrangian $\tildeN^\ast_k$ after perturbation since this would make
$\CT(\tildeN^\ast_k)$ a Lagrangian submanifold of the resolved
conifold. Unfortunately there is an obstruction to such perturbation
following from a theorem of Gromov-Oh: in a cotangent bundle every exact Lagrangian submanifold intersects
the zero section \cite{Gr,ALP,Oh1}.

Thus we have to settle for an ordinary isotopy instead of a
symplectic one. Even though $\CT(\tildeN^\ast_k)$ will no longer be Lagrangian
in $\hatcalC$, it will be good enough for the purposes of open
Gromov-Witten theory. Specifically, it will be Lagrangian with respect to a different uniformly tame symplectic form on $\calC$ (tame Lagrangian, see \refD{5}).

As before we identify $T^*S^3\simeq TS^3$. To separate $N^\ast_k$ from the zero section we simply move it within each
fiber in the direction tangent to the knot. Recall that
$N^\ast_k= \{ (k(t),p)\in T^\ast S^3\mid t\in S^1,p\cdot\dotk(t)=0\}$.

\begin{definition}\label{D:3}
The perturbed conormal bundle is
\be\label{e3.1}
  N^\ast_{k,\ve}:=\{(k(t),p+\ve\dotk(t))\mid t\in S^1,
  \quad p\cdot\dotk(t)=0\} 
\ee
\end{definition}
Since $\dotk(t)\cdot k(t)=0$one has $N^\ast_{k,\ve}\subset T^\ast S^3$ for all $\ve\geq0$.
And $p\cdot\dotk(t)=0$ implies $|p+\ve\dotk(t)|^2=|p|^2+\ve^2|\dotk(t)|^2=|p|^2+\ve^2\geq\ve^2>0$
so $N^\ast_{k,\ve}$ is indeed disjoint from the zero section for any $\ve >0$.

Since $\CT(N^\ast_{k,\ve})=\pi_2^{-1}\circ F(N^\ast_{k,\ve})$ the proofs in this section split into two parts: first we prove that $F(N^\ast_{k,\ve})$ is a tame Lagrangian in $\calC$ and then that $\pi_2^{-1}$ preserves this property. \refL{7} is also used to prove bounded geometry of $\hatcalC$ in the next section.
It is convenient to represent $N^\ast_{k,\ve}$  as the image of
$N^\ast_k$ under an ambient isotopy in $\R^4\times\R^4\supset T^\ast S^3$.
To this end, let $\xi$ be a smooth vector field in $\R^4$ with compact support
satisfying $\xi(k(t))=\dotk(t)$.
\begin{definition}\label{D:4}
Let 
\bee
  \begin{aligned}
  \R^4\times\R^4   & \overset{\Phi_\ve}{\lra} \calC_a\\
  (x,p)        &\longmapsto (x,p+\ve\xi(x)).
 \end{aligned}
\eee
$\Phi_\ve$ is an isotopy since
$\Phi^{-1}_\ve(x,p)=(x,p-\ve\xi(x))$ and by
construction of $\xi$, $\Phi_\ve(N^\ast_k)=N^\ast_{k,\ve}$.
\end{definition}

\begin{figure}\label{F:4}
\hskip2in\epsfig{file=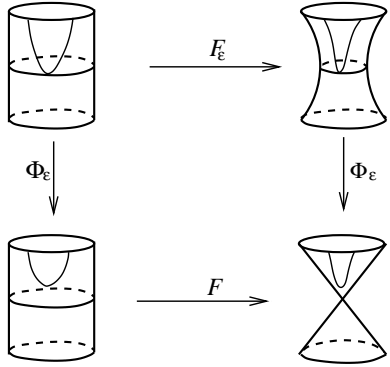,width=2truein}
\caption{Disconnecting the conormal bundle from the zero section.}
\end{figure}

Let us now look at the image of $N^\ast_{k,\ve}$ under the contraction
$F(x,p)=(|p|x,p)$, the first half of the conifold transition:
\bee
  F(N^\ast_{k,\ve})=
  \left\{\left(x\sqrt{|p|^2+\ve^2},p+\ve\dotk(t)\right)
  \mid t\in S^1, p\cdot\dotk(t)=0\right\}.
  \eee
Recall the regularized contraction $F_\ve$ from \refS{1}: $F_\ve(x,p):=\left(x\sqrt{|p|^2+\ve^2},p\right)$.
The following diagram of maps commutes (see Fig.\,4):
$$ 
\begin{diagram}
N^\ast_k & \rTo^{F_\ve} & \R^4\times\R^4 \\
\dTo^{\Phi_\ve}&              & \dTo_{\Phi_\ve}\\ 
\R^4\times\R^4 & \rTo^F       & \R^4\times\R^4.
\end{diagram} 
$$
Now we want to investigate in what sense
$F(N^\ast_{k,\ve})=F(\Phi_\ve(N^\ast_k))=\Phi_\ve(F_\ve(N^\ast_k))$
is close to being Lagrangian in $\calC$.

\begin{lemma}\label{L:3}
$F_\ve(N^\ast_k)$ is isotropic in $\R^4\times\R^4$, i.e.,
$\omega=\Sigma_i dx_i\wedge dp_i$ vanishes on it.
\end{lemma}

\begin{proof}
From the parametrization of $N^\ast_k$ we have
for $F_\ve(N^\ast_k)$:
\bee
  \begin{aligned}
  S^1\times\R^2    & \overset{F_\ve\circ K}{\longmapsto} \R^4\times\R^4\\
  (t,\alpha,\beta) & \longmapsto
    (k(t)\sqrt{\alpha^2+\beta^2+\ve^2},\alpha p^1(t)+\beta p^2(t)).
  \end{aligned}\eee
Thus the tangent bundle is spanned by
\bee
  \begin{aligned}
  (F_\ve\circ K)_\ast \partial_f
    & =(\dotk(t)\sqrt{\alpha^2+\beta^2+\ve^2},
       \alpha\dotp^1(t)+\beta\dotp^2(t)) \\
  (F_\ve\circ K)_\ast \partial_\alpha
    & =(k(t)\alpha(\alpha^2+\beta^2+\ve^2)^{-1/2}, p^1(t)) \\
  (F_\ve\circ K)_\ast \partial_\beta
    &=(k(t)\beta(\alpha^2+\beta^2+\ve^2)^{-1/2}, p^2(t))\\
  \end{aligned}
\eee
Let $J$ be the standard complex structure on $\R^4\times\R^4\simeq\C^4$, then
\bee
  \begin{aligned}
   J (F_\ve\circ K)_\ast\partial_t
    &=(-\alpha\dotp^1(t)-\beta\dotp^2(t),
      (\alpha^2+\beta^2+\ve^2)^{1/2}\dotk(t))\\
   J (F_\ve\circ K)_\ast\partial_\alpha
    &=(-p^1(t),\alpha(\alpha^2+\beta^2+\ve^2)^{-1/2}k(t))\\
   J (F_\ve\circ K)_\ast\partial_\beta
    &=(-p^2(t),\beta(\alpha^2+\beta^2+\ve^2)^{-1/2}k(t)).
  \end{aligned}\eee
Since $k(t)$, $\dotk(t)$, $p^1(t)$, $p^2(t)$ are pairwise orthogonal for
each $t$ one can see by inspection that
$ J (T F_\ve(N^\ast_k))\perp T F_\ve(N^\ast_k)$, which is
equivalent to isotropy.
\end{proof}

However, $F(N^\ast_{k,\ve})=\Phi_\ve(F_\ve(N^\ast_k))$ will no longer be
Lagrangian in $\calC$. It does satisfy a weaker property that we now introduce \cite{MS,S}.
\begin{definition}\label{D:5}
Let $(M, J ,g)$ be an almost K\"ahler manifold.
A symplectic form $\tildeo$ is called uniformly tame if there exists a constant $C\geq 1$
such that for any vector field $X$ on $M$:
\be\label{e3.2}
C^{-1}g(X,X)\leq \tildeo(X, JX)\leq C g(X,X) 
\ee
A submanifold $L\hra M$ is (uniformly) tame isotropic if there
is a uniformly tame $\tildeo$ defined in its neighborhood such that
$\tildeo\mid_{TL}=0$.
If in addition $\dim L=\frac{1}{2}\dim M$ then $L$ is called (uniformly)
tame Lagrangian.
\end{definition}
Note that the K\"ahler form $\omega(X,Y)=g(JX,Y)$ is obviously
uniformly tame with $C=1$ so this is a generalization of the Lagrangian
condition. Also, tame Lagrangian implies totally real, i.e.,
$J(TL)\bigcap TL=\{0\}$.
Indeed, if $X\in J(TL)\bigcap TL$, then $JX\in TL$
since $J^2=-I$ and
$|X|^2=g(X,X)\leq C\tildeo(X,JX)=0$ so $X=0$.
Most importantly for us, the property of being tame Lagrangian or
isotropic is preserved under biholomorphisms. Namely, if
$\Phi:(M,J,g)\to(\tildeM,\tildeJ,\tildeg)$
is a biholomorphism and $L\hra M$ is $\tildeo$-Lagrangian
then $\Phi(L)$ is $\Phi_\ast\tildeo:=\Phi^{-1\ast}\tildeo\ $-Lagrangian.

So far we have not imposed any restrictions on the perturbation parameter
$\ve$. We will do so now to ensure uniform tameness.
Let $\xi$ be the perturbation vector field from \refL{3} and $D\xi$ be its
Jacobian matrix. Set $\|D\xi\|:=\sup_{x\in\R^4}\|D\xi(x)\|$

\begin{lemma}\label{L:4}
If $\ve<\frac{1}{\|D\xi\|}$,
then $\Phi_{\ve\ast}\omega=\Phi^{-1\ast}_\ve\omega$ is uniformly tame in
$\R^4\times\R^4\simeq\C^4$.
\end{lemma}

\begin{proof}
\bee\label{e3.8}
  \begin{aligned}
  \Phi_{\ve\ast}\omega=\Phi^{-1\ast}_\ve \omega
  &=\Sum_i\Phi^{-1\ast}_\ve(dx_i\wedge dp_i)
  =\Sum_i dx_i\wedge d(p_i-\ve\xi_i(x))\\
  &=\Sum_i dx_i\wedge dp_i -\ve\Sum_{i,j}
    \frac{\partial\xi_i}{\partial x_j}dx_i\wedge dx_j
  \end{aligned}
  \eee
Since $J\partial_{x_i}=\partial_{p_i}$
and $J\partial_{p_i}=-\partial_{x_i}$, we get
\bee\label{e3.9}
\begin{aligned}
  \Phi_{\ve\ast} \omega(X,JX) &=\Sum_i(dx_i(X)^2+dp_i(X)^2)-\ve\Sum_{i,j}
  \left( \frac{\partial\xi_i}{\partial x_j}
  -\frac{\partial\xi_i}{\partial x_j} \right)
   dp_i(X) dx_j(X)\\
   &=|X|^2-\ve \left((D\xi-D\xi^T) dp(X), dx(X) \right)
   \end{aligned}\eee
The absolute value of the second term is bounded by
\bee\label{e3.10}
  2\ve\|D\xi\| |dp(X)| |dx(X)| \leq 2\ve\|D\xi\|
    \frac{|dp(X)|^2 + |dx(X)|^2}{2}=\ve\|D\xi\||X|^2
  \eee
Therefore
$(1-\ve\|D\xi\|)|X|^2\leq\Phi_{\ve\ast}
  \omega(X,JX)\leq(1+\ve\|D\xi\|)|X|^2$
and if $\ve<\frac{1}{\|D\xi\|}$, $\Phi_{\ve\ast}\omega$ is uniformly tame
with $C=(1-\ve\|D\xi\|)^{-1}$
\end{proof}

\begin{corollary}\label{C:1}
$F(N^\ast_{k,\ve})$ is a tame Lagrangian in
$\calC$ for $\ve<\frac{1}{\|D\xi\|}$.
\end{corollary}

\begin{proof}
Let $\tildeo=\Phi_{\ve\ast}\omega$. Since $\Phi_\ve$ is a diffeomorphism
and $\Phi_\ve(F_\ve(N^\ast_k))=F(N^\ast_{k,\ve})$ vectors $\Phi_{\ve\ast}X$
for $X\in TF_\ve(N^\ast_k)$ span $TF(N^\ast_{k,\ve})$. But
\bee
  \tildeo(\Phi_{\ve\ast}X,\Phi_{\ve\ast}Y)
  =\Phi^{-1\ast}_\ve\omega(\Phi_{\ve\ast}X,\Phi_{\ve\ast}Y)
  =\omega(X,Y)=0 \eee
since $F_\ve(N^\ast_k)$ is isotropic by \refL{3}.
\end{proof}

As $N^\ast_{k,\ve}$ is disjoint from the zero section of
$T^\ast S^3$ its contraction $F(N^\ast_{k,\ve})$ now avoids the singularity
of the conifold at the origin.
In fact, if $(x,p)\in F(N^\ast_{k,\ve})$ then $|x|^2+|p|^2\geq 2\ve^2>0$.
Now the second half of the conifold transition \eqref{e1.4} constitutes lifts $F(N^\ast_{k,\ve})$ to $\hatcalC$. Recall from \eqref{e1.01} that we changed the coordinates in $\C^4$ from $z=x+ip$ to $w=\sqrt{2}(Uz)$, where $U$ is a unitary matrix. Therefore, along $F(N^\ast_{k,\ve})$ we have $|w|\geq\sqrt{2}\cdot \sqrt{2}\ve=2\ve$ and it is separated from $0$. Since $\pi^{-1}_2:\calC\setminus\{0\}\to\hatcalC\setminus 0(\hatcalC)$
is a biholomorphism there is no need to take closure in the conifold transition \eqref{e1.04} and we simply have 
$$
\CT(N^\ast_{k,\ve})=\pi^{-1}_2(F(N^\ast_{k,\ve})).
$$
To establish that $\CT(N^\ast_{k,\ve})$ is a tame Lagrangian in $\hatcalC$ we need certain properties of $\pi^{-1}_2$.

\begin{definition}[\cite{ALN}]\label{D:6}
A map $\Phi:(M,g)\to(\tildeM,\tildeg)$ between Riemannian manifolds is
(uniformly) bi-Lipschitz if there exists a constant $C\geq1$ such that 
\be\label{e3.3}
  C^{-1}g\leq \Phi^\ast\tildeg\leq Cg. 
\ee
\end{definition}
The map $\Phi$ does not have to be a diffeomorphism, e.g., any isometric immersion
would satisfy this condition with $C=1$.
If $\Phi$ is a diffeomorphism then being uniformly bi-Lipschitz is
equivalent to the norms of $\Phi_\ast$ and $\Phi^{-1}_\ast$ being bounded by a constant,
i.e., $|\Phi_\ast X|_\tildeg\leq C|X|_g$
and $|\Phi^{-1}_\ast Y|_g\leq C|Y|_\tildeg$.

\begin{lemma}\label{L:5}
Let $\Phi:(M,J,g)\to(\tildeM,\tildeJ,\tildeg)$ be a bi-Lipschitz
biholomorphism between two almost K\"ahler manifolds.
If $L\hra M$ is a tame Lagrangian then so is $\tildeL:=\Phi(L)$.
\end{lemma}

\begin{proof}
Let $\omega'$ be the corresponding uniformly tame symplectic form in the
neighborhood of $L$ in $M$, i.e., $\omega'\mid_{TL}=0$,
$C^{-1}_1 g(X,X)\leq \omega'(X,JX)\leq C_1g(X,X)$.
Consider $\Phi_\ast\omega'=\Phi^{-1\ast}\omega'=:\tildeo'$ in the
neighborhood of $\tildeL$. Since $\Phi$ is a biholomorphism:
\bee\label{e3.13}
  \tildeo'(\Phi_\ast X,J\Phi_\ast X)=\tildeo'(\Phi_\ast X,\Phi_\ast JX)
  =\Phi^{-1\ast}\omega'(\Phi_\ast X,\Phi_\ast JX)=\omega'(X,JX).
  \eee
Since $\Phi$ is bi-Lipschitz
\bee\label{e3.14}
  C^{-1}_2 g(X,X)\leq\tildeg(\Phi_\ast X,\Phi_\ast X)\leq C_2 g(X,X). \eee
Therefore, combining the inequalities
\bee\label{e3.15}
  (C_1C_2)^{-1}\tildeg(\Phi_\ast X,\Phi_\ast X)
  \leq\tildeo'(\Phi_\ast X,J\Phi_\ast X)
  \leq C_1C_2\,\tildeg(\Phi_\ast X,\Phi_\ast X),
  \eee
and $\tildeo'$ is uniformly tame. Also if $X,Y\in TL$
\bee
  \tildeo'(\Phi_\ast X,\Phi_\ast Y)
  =\Phi^{-1\ast}\omega'(\Phi_\ast X,\Phi_\ast Y)=\omega'(X,Y)=0
  \eee
and $\tildeo'\mid_{T\tildeL}=0$
\end{proof}

In view of \refL{5} to prove that $\CT(N^\ast_{k,\ve})$ is a tame Lagrangian 
we have to show that $\pi^{-1}_2$ is bi-Lipschitz away
from the singularity in $\calC$. Recall that the metric on $\hatcalC$ is 
$\hatg=\pi^\ast_1 g_{FS}+\pi^\ast_2 g_{st}$,
where $g_{FS}$ is the Fubinin-Study metric on $\CP^1$ and $g_{st}$
is the standard metric on $\C^4$.
Hence
$$
\pi^{-1\ast}_2\hatg=(\pi_1\circ\pi^{-1}_2)^\ast g_{FS}+g_{st}\geq g_{st}.
$$
To establish the inverse inequality let us introduce convenient notation.
\begin{definition}\label{D:7}
If $\alpha,\beta$ are 1-forms then
$\alpha\odot\beta:=\alpha\otimes\beta+\beta\otimes\alpha$
defines their symmetric product.
\end{definition}
Then for instance
$g_{FS}=\frac{1}{2}\frac{dz\odot\bardz }{(1+|z|^2)^2}$,
where $z=\frac{u}{v}$ and $[u:v]$ are homogeneous coordinates on $\CP^1$.
One also has a Cauchy inequality:
\be\label{e3.4}
|\alpha\odot\barb|\leq\frac{1}{2}(\alpha\odot\bara+\beta\odot\barb)
\ee

\begin{lemma}\label{L:6} With the above notation
\be\label{e3.5}
  (\pi_1\circ\pi^{-1}_2)^\ast g_{FS}
  \leq \frac{2}{|w|^2} g_{st}\ \text{{\rm for}}\ w\in\calC\setminus\{0\}.
\ee
\end{lemma}

\begin{proof}
If
$([u:v],w_1,w_2,w_3,w_4)\in\hatcalC$ then
$\begin{vmatrix}w_1&w_2\\u&v\end{vmatrix}
=\begin{vmatrix}w_3&w_4\\u&v\end{vmatrix}=0$.
Therefore $z=\frac{u}{v}=\frac{w_1}{w_2}$ if $w_2\not=0$
and $\frac{w_3}{w_4}$ if $w_4\not=0$.
Assume for now that $w_2\not=0$ and $(\pi_1\circ\pi^{-1}_2)(w)=[w_1:w_2]$.
Then
\bee
  \begin{aligned}
  (\pi_1\circ\pi^{-1}_2)^\ast
  &  g_{FS}
     =\frac{1}{2}\frac{d\frac{w_1}{w_2}\odot\overline{d\frac{w_1}{w_2}}}
       {(1+|\frac{w_1}{w_2}|^2)^2}\\
  &=\frac{1}{2}
    \frac{ (w_2 dw_1-w_1dw_2)\odot\overline{(w_2dw_1-w_1dw_2)}}
         { (|w_1|^2+|w_2|^2)^2 }\\
  &=\frac{1}{2}
    \frac{ |w_2|^2dw_1\odot\overline{dw_1}+|w_1|^2dw_2\odot\overline{dw_2}
             -w_1\overline{w_2}dw_2\odot\overline{dw_1}
             -w_2\overline{w_1}dw_1\odot\overline{dw_2} }
         {(|w_1|^2+|w_2|^2)^2}.
   \end{aligned} \eee
By the Cauchy inequality (\ref{e3.4}):
\bee
  |w_1\overline{w_2} dw_2\odot \overline{dw_1}|
    =|w_1dw_2\odot\overline{w_2dw_1}| \leq\frac{1}{2}
  \left( |w_1|^2 dw_2\odot \overline{dw_2}
   + |w_2|^2 dw_1\odot\overline{dw_1} \right) \eee
and the same estimate holds for the second cross-term.
Hence
\bee
  \begin{aligned}
  (\pi_1\circ\pi^{-1}_2)^\ast g_{FS}
  &\leq \frac{|w_2|^2dw_1\odot\overline{dw_1}+|w_1|^2dw_2\odot\overline{dw_2}}
             {(|w_1|^2+|w_2|^2)^2} \\
  &\leq \frac{dw_1\odot \overline{dw_1}+dw_2\odot\overline{dw_2}}
             {(|w_1|^2+|w_2|^2)},
  \end{aligned}\eee
or
$$
(|w_1|^2+|w_2|^2)(\pi_1\circ\pi^{-1}_2)^\ast g_{FS}\leq dw_1\odot\overline{dw_1}+dw_2\odot\overline{dw_2}.
$$
Analogously, if $w_4\not=0$
\bee
  (|w_3|^2+|w_4|^2)(\pi_1\circ\pi^{-1}_2)^\ast g_{FS}
  \leq dw_3\odot \overline{dw_3}+dw_4\odot\overline{dw_4} \eee
Adding together the last two inequalities and taking into account that
$g_{st}=\frac{1}{2}\Sum^4_{i=1}dw_i\odot\overline{dw_i}$, one gets
\bee
  |w|^2(\pi_1\circ\pi^{-1}_2)^\ast g_{FS}\leq 2 g_{st}. \eee
Although we assumed $w_2\not=0$, $w_4\not= 0$ in the process,
the final inequality holds by continuity for any $w\in\calC$ and for
$w\not=0$: $(\pi_1\circ\pi^{-1}_2)^\ast g_{FS}\leq\frac{2}{|w|^2}g_{st}$.
\end{proof}

\begin{corollary}\label{C:2}
$\pi^{-1}_2$ is uniformly bi-Lipschitz on any $\calC\setminus B_R(0)$
with $C=1+\frac{2}{R^2}$.
\end{corollary}

\begin{proof}
$g_{st}\leq g_{st}+(\pi_1\circ\pi^{-1}_2)^\ast g_{FS}=\pi^{-1\ast}_2\hatg
\leq(1+\frac{2}{R^2})g_{st}$  by \refL{6}.
\end{proof}

Now we are ready to state the main result of this section:

\begin{theorem}\label{T:2}
Let $k:S^1\to S^3\hra\R^4$ be a knot and
$N^\ast_{k,\ve}$ be its perturbed conormal bundle. If $\ve<1$ then
$\CT(N^\ast_{k,\ve})\hra\hatcalC$ is a tame Lagrangian
in the resolved conifold. Moreover, the form $\tildeo_\ve:=\pi^{-1}_{2\ast}\Phi_{\ve\ast}\omega$ 
that makes it Lagrangian is exact on $\hatcalC\setminus 0(\hatcalC)$.
\end{theorem}

\begin{proof}
By \refC{1} $F(N^\ast_{k,\ve})$ is a tame Lagrangian in $\calC$ if
$\ve<\frac{1}{\|D\xi\|}$, where $\xi$ extends $\dotk(t)$ from $k$ to $\R^4$.
Since $|\dotk(t)|=1$, the extension can be carried out so that $\|D\xi\|=1$.
$F(N^\ast_{k,\ve})\subset\calC\setminus B_{2\ve}(0)$ and $\pi^{-1}_2$
is a bi-Lipschitz  biholomorphism on any $\calC\setminus B_R(0)$ by
\refC{2}. Since $\CT(N^\ast_{k,\ve})=\pi^{-1}_2 (F(N^\ast_{k,\ve}) )$
the latter is a tame Lagrangian by \refL{5}. Since $\tildeo_\ve$ is a pushforward and 
therefore a pullback of an exact form on $T^\ast S^3$ it is exact.
\end{proof}

\section{Geometry of the conifold}\label{S:4}

The ultimate goal of constructing Lagrangian or totally real submanifolds
of the resolved conifold is to consider the moduli of holomorphic curves
ending on them and to define open Gromov-Witten invariants. At the very least, one
needs these moduli spaces to be compact. Since neither the resolved conifold itself nor the conormal bundles and their conifold transitions are compact certain uniform bounds are required to ensure compactness of the moduli. They are known as bounded geometry \cite{S,ALP} and generalize geometric properties of closed Riemannian manifolds and their submanifolds.

\begin{definition}\label{D:8}
A Riemannian manifold $(M,g)$ has bounded
geometry (or is geometrically bounded) if its sectional curvature is bounded from above
$\sec(X,Y)\leq K<\infty$ and its injectivity radius is bounded from below $i(M)\geq r_0>0$.
\end{definition}

The main result of this section (\refT{3}) claims that the resolved
conifold $\hatcalC$ is geometrically bounded. We use it to prove compactness
of the moduli of holomorphic curves ending on $\CT(N^\ast_{k,\ve})$ in the
next section.

To obtain estimates on curvature it is convenient to use second fundamental
forms. Recall the definition \cite{PdC}:

\begin{definition}\label{D:9}
Let $L\hra(M,g)$ be a smooth submanifold and
$X,Y\in\Gamma(TL)$ be vector fields on it. Then
\bee
  \II^{L/M}(X,Y):=\pr_{T^\perp L}(\nabla_X Y), 
  \eee
where $\nabla$ is the Riemannian connection on $M$ and
$\pr_{T^\perp L}$ is the orthogonal projection to the normal bundle $T^\perp L$ of $L$ in $TM$ is called the second fundamental form of $L$ in $M$.
\end{definition} 

If $L\hra Q\hra (M,g)$, it follows from linear
algebra that $\II^{L/M}=\II^{L/Q}+\II^{Q/M}$
and the terms on the right are orthogonal to each other.
In particular, 
\be\label{e4.1}
|\II^{L/Q}| \leq |\II^{L/M}|.
\ee
When there is no confusion we drop the ambient manifold from notation and
write simply $\II^L$. The norm $\|\II^L\|$ is the smallest number $C$
such that $|\II^L(X,Y)\mid_g \leq C|X|_g|Y|_g$. Our interest in the second 
fundamental forms is explained by the Gauss equation \cite{PdC}. If $X,Y\in\Gamma(TL)$ and
$\sec^L(X,Y)$, $\sec^M(X,Y)$ denote the sectional curvatures in
$L$ and $M$ respectively:
\be\label{e4.2}
  \sec^L(X,Y)=\sec^M(X,Y)+g(\II^L(X,X),\II^L(Y,Y))
  -|\II^L(X,Y)|^2_g. 
\ee
Thus a bound on ambient curvature and second fundamental form yields one
on the curvature of a submanifold. To consider behavior of second
fundamental forms under smooth maps we need the following notion.

\begin{definition}[\cite{EL}]\label{D:10}
Let $\Phi:(M,g)\to(\tildeM,\tildeg)$ be a smooth map between two
Riemannian manifolds. Let $\nabla,\tildenabla$ be the respective
Riemannian connections and $X,Y\in\Gamma(TM)$.
The covariant Hessian  $\nabla^2\Phi$ ('second fundamental form of a map' in \cite{EL} ) is by definition
\be\label{e4.3}
  \nabla^2\Phi(X,Y):=\tildenabla_{\Phi_\ast X}(\Phi_\ast Y)-\Phi_\ast(\nabla_XY). 
\ee
$\nabla^2\Phi\in\Gamma((T^\ast M)^{\otimes 2}\otimes\Phi^\ast(T\tildeM))$
and it is straightforward to check that it is symmetric and tensorial.
$\|\nabla^2\Phi\|$ is the smallest number $C$ such that
$|\nabla^2\Phi(X,Y)|_\tildeg \leq C|X|_g|Y|_g$.
\end{definition}

\begin{lemma}\label{L:7}
Let $X,Y\in\Gamma(TL)$ and
$$
\begin{diagram}
(M,g)  & \rTo^\Phi &  (\tildeM,\tildeg)  \\
\uInto &           &  \uInto        \\
 L     & \rTo^\Phi &  \tilde{L}    
\end{diagram}
$$
Then 
\be\label{e4.5}
 \II^{\tilde{L}}(\Phi_\ast X,\Phi_\ast Y)=\pr_{T^\perp L}
 (\Phi_\ast \II^L(X,Y) + \nabla^2\Phi(X,Y))
\ee
\end{lemma}  
\begin{proof}
\bee
  \begin{aligned}
  \II^{\tildeL}(\Phi_\ast X,\Phi_\ast Y)
  &=\pr_{T^\perp L} (\tildenabla_{\Phi_\ast X} \Phi_\ast Y) \\
  &=\pr_{T^\perp L}
    ( (\tildenabla_{\Phi_\ast X}\Phi_\ast Y -\Phi_\ast(\nabla_XY))
    +\Phi_\ast(\nabla_X Y)) \\
  &=\pr_{T^\perp L}
    (\nabla^2\Phi(X,Y)+\Phi_\ast(\pr_{T^\perp L}(\nabla_X Y)
    +\pr_{TL} (\nabla_X Y)))\\
  &=\pr_{T^\perp L}
    (\nabla^2\Phi(X,Y)+\Phi_\ast \II^L(X,Y))
    +\pr_{T^\perp L}(\Phi_\ast\pr_{TL}(\nabla_X Y)) \\
  &=\pr_{T^\perp L}
    (\Phi_\ast \II^L(X,Y)+\nabla^2\Phi(X,Y))
  \end{aligned}\eee
since $\Phi_\ast(TL)\subset T\tilde{L}$.
\end{proof}

\begin{corollary}\label{C:3}
If $\Phi:(M,g)\to(\tildeM,\tildeg)$ is a bi-Lipschitz map with a bounded
covariant Hessian then images of submanifolds with bounded second
fundamental forms have bounded second fundamental forms.
\end{corollary}

\begin{proof}
Recall that by definition of a bi-Lipschitz map
\bee
  C^{-1}g\leq \Phi^\ast\tildeg\leq Cg \eee
and in particular
$|\Phi_\ast X|^2_\tildeg =\tildeg(\Phi_\ast X,\Phi_\ast X)
 =\Phi^\ast \tildeg(X,X)\leq C|X|^2_g$
so $\|\Phi_\ast\|\leq C^{1/2}$. By \refL{7}
\bee
  \begin{aligned}
  |\II^{\tilde{L}}(\Phi_\ast X,\Phi_\ast Y)|_g
  &\leq|\nabla^2\Phi(X,Y)|_\tildeg+|\Phi_\ast \II^L(X,Y)\mid_\tildeg \\
  &\leq\|\nabla^2\Phi\| |X|_g|Y|_g+\|\Phi_\ast||\|\II^L\| |X|_g |Y|_g\\
  &\leq(\|\nabla^2\Phi\|+C^{1/2}\|\II^L\|)|X|_g |Y|_g \\
  &\leq C(\|\nabla^2\Phi\|+C^{1/2}\|\II^L\|)
    |\Phi_\ast X|_\tildeg |\Phi_\ast Y|_\tildeg
  \end{aligned}\eee
Thus $\|\II^{\tilde{L}}\| \leq C(\|\nabla^2\Phi\| + C^{1/2} \|\II^L\|)$
\end{proof}

Note that if $\Phi$ is a bi-Lipschitz embedding, then $L=M$,
$\II^{L/M}=\II^{M/M}=0$ and $\|\II^{\tilde{L}}\|\leq C\|\nabla^2\Phi\|$.
In other words, second fundamental forms of embeddings are controlled by
the bi-Lipschitz constants and covariant Hessians.

To obtain estimates on covariant Hessians, we need a coordinate
representation. For convenience we use the following notation:

\bee
  \begin{aligned}
  \partial\Phi(X)
   &=\frac{\partial\Phi}{\partial x_i}X^i, \quad
  \partial^2\Phi(X,Y)=\frac{\partial^2\Phi}{\partial x_i\partial x_j}X^iY^i\\
  \Gamma(X,Y)
   &=\Gamma^k_{i_j}X^i Y^j, \quad
   \tildeG(\tildeX,\tildeY)=\tildeG^k_{i_j}\tildeX^i\tildeY^j,
  \end{aligned}\eee
where $\Gamma^k_{ij}$, $\tildeG^k_{ij}$ are Christoffel symbols for
$\nabla$, $\tildenabla$ respectively.

\begin{lemma}\label{L:8}
Let $\Phi:(M,g)\to(\tildeM,\tildeg)$. Then in local coordinates
\be\label{e4.6}
  \nabla^2\Phi(X,Y)=\partial^2\Phi(X,Y)
  +\tildeG(\partial\Phi(X),\partial\Phi(Y))
  -\partial\Phi(\Gamma(X,Y)).
\ee
\end{lemma}

\begin{proof}
Since $\nabla^2\Phi$ is tensorial, we can ignore expressions containing
derivatives of $X$, $Y$ in the calculation:
\bee
  \begin{aligned}
  \nabla^2\Phi(X,Y)
  &=\tildenabla_{X^i\Phi_\ast\partial_{x_i}}(Y^j\Phi_\ast\partial_{x_j})
    -\Phi_\ast (\nabla_{X^i\partial_{x_i}} Y^j\partial_{x_j}) \\
  &=X^iY^j\tildenabla_{\Phi_\ast\partial_{x_i}} \Phi_\ast \partial_{x_j}
    -X^iY^j\Phi_\ast (\nabla_{\partial_{x_i}} \partial_{x_j}) \\
  &=X^iY^j (\tildenabla_{\Phi_\ast\partial_{x_i}}
    \left(\frac{\partial\Phi^\alpha}{\partial_{x_j}}\partial y_\alpha\right)
    -\Phi_\ast (\Gamma^k_{ij}\partial_{x_k})) \\
  &=X^iY^j
  \left(\frac{\partial^2\Phi^\alpha}{\partial x_i\partial x_j}\partial_{y_\alpha}
   +\frac{\partial\Phi^\alpha}{\partial x_j}
   \tildenabla_{\frac{\partial\Phi^\beta}{\partial x_i} \partial_{y_\beta}}
   \partial_{y_\alpha}-\frac{\partial\Phi^\alpha}{\partial x_k}
   \Gamma^k_{ij} \partial_{y_\alpha} \right)  \\
  &=X^iY^j(
  \left(\frac{\partial^2\Phi^\alpha}{\partial x_i\partial x_j}\partial_{y_\alpha}
   +\frac{\partial\Phi^\alpha}{\partial x_j}
   \frac{\partial\Phi^\beta}{\partial x_i} \tildeG^\gamma_{\beta\alpha}
   \partial_{y_\gamma} -\frac{\partial\Phi^\alpha}{\partial x_k}
   \Gamma^k_{ij} \partial_{y_\alpha} \right)   \\
  &=\left( \frac{\partial^2\Phi^\alpha}{\partial x_i\partial x_j}
       X^iY^j +\tildeG^\alpha_{\beta\gamma}
       \left(\frac{\partial\Phi^\beta}{\partial x_i}X^i \right)
       \left(\frac{\partial\Phi^\gamma}{\partial x_j}Y^j\right)
       -\frac{\partial\Phi^\alpha}{\partial x_k}
       \left(\Gamma^k_{ij} X^iY^j\right)
     \right)\partial_{y_\alpha} \\
  &=\partial^2\Phi(X,Y)+\tildeG(\partial\Phi(X),\partial\Phi(Y))
    -\partial\Phi(\Gamma(X,Y))
  \end{aligned}\eee
\end{proof}

\noindent Note that for a map between two flat spaces $\nabla^2\Phi$ turns into the
usual Hessian $\partial^2\Phi$.

Now we turn to the geometry of the conifold
\bee
  \calC=\{w\in\C^4\mid w_1w_4=w_2w_3\}.
\eee
We start by finding convenient parametrizations for
$\calC\setminus\{0\}$.

\begin{lemma}\label{L:9}
Let $w\in\calC\setminus\{0\}$.
Then there exist $\xi,\eta\in\C$ and
$z\in\overline{\D}:=\{z\in\C\mid |z|\leq 1\}$ such that
$w=(\xi,z\xi,\eta,z\eta)$ or $w=(z\xi,\xi,z\eta,\eta)$.
Moreover, 
\bee
  \begin{aligned}
  \C\times(\C^2\setminus\{0\})
    &\overset{\Phi}{\lra}\C^4\\
  (z,\xi,\eta)
    &\longmapsto (\xi,z\xi,\eta,z\eta)
  \end{aligned}
\eee
is an embedding.
\end{lemma}
\begin{proof}
Since $w\not=0$ at least one of $w_i$ is non-zero, let $w_1\not=0$.
If $|w_2|>|w_1|$ then $w_2\not=0$ and we can set $\xi=w_2$,
$z=\frac{w_1}{w_2}$ with $|z|<1$.
Otherwise set $\xi=w_1$, $z=\frac{w_2}{w_1}$ with $|z|\leq 1$.
Since the two cases are analogous, let us consider just the latter one.
In this case $w_3=0$ implies $w_4=0$ since $w_1w_4=w_2w_3$
and $w_1\not=0$.
Therefore $w=(\xi,z\xi,\eta,z\eta)$ with $\xi,z$ as above
and $\eta=0$. If $w_3\not=0$, then $\frac{w_4}{w_3}=\frac{w_2}{w_1}=z$
and $w=(\xi,z\xi,\eta,z\eta)$ with $\eta=w_3$.
The second possibility of $w=(z\xi,\xi,z\eta,\eta)$ arises when
$|w_2|>|w_1|$.

The Jacobian of $\Phi$ is
$\begin{pmatrix}0&\xi&0&\eta\\1&z&0&0\\0&0&1&z\end{pmatrix}$
and it obviously has full rank unless $\xi=\eta=0$ so $\Phi$ is an immersion.
Also, if $(\xi,z\xi,\eta,z\eta)=(\xi',z'\xi',\eta',z'\eta')$
then $\xi=\xi'$, $\eta=\eta'$ and one of them, say, $\xi=\xi'\not=0$.
But then $z'=\frac{z'\xi'}{\xi'}=\frac{z\xi}{\xi}=z$ and $\Phi$
is an embedding.
\end{proof}      

\refL{9} implies that $\calC\setminus\{0\}$ can be covered by two
charts, each of them diffeomorphic (in fact, biholomorphic) to
$2\D\times(\C^2\setminus\{0\})$.
Since the corresponding parametrizations are the same up to permutation of
coordinates we may consider just one of them. First of all we want to
describe the induced metric on $\calC$ in terms of $z$, $\xi$, $\eta$.

\begin{lemma}\label{L:10}
Let 
\bee
  \begin{aligned}
  2\D\times(\C^2\setminus\{0\})
    &\overset{\Phi}{\lra}\C^4\\
  (z,\xi,\eta)
    &\longmapsto (\xi,z\xi,\eta,z\eta)
  \end{aligned}
\eee
and
\bee
  g:=\frac{1}{2}\left(
  (|\xi|^2 +|\eta|^2)dz\odot\bardz  +d\xi\odot\bardxi
  +dy\odot\overline{dy}\right)
\eee
be a metric on $2\D\times(\C^2\setminus\{0\})$.
Then $\Phi$ is uniformly bi-Lipschitz.
\end{lemma}

\begin{proof}
\bee
 \begin{aligned}
 \Phi^\ast g_{st}
 &=\frac{1}{2}
   \left(d\xi\odot\bardxi+d(z\xi)\odot\overline{d(z\xi)}
   +d\eta\odot\bardeta+d(z\eta)\odot\overline{d(z\eta)}\right)\\
 &=\frac{1}{2}
   \left( d\xi\odot\bardxi+|\xi|^2dz\odot\bardz 
   +|z|^2d\xi\odot\bardxi+|\eta|^2dz \odot\bardz 
      +|z|^2d\eta\odot\bardeta \right.\\
 & \left.\qquad
   +\xi\barz dz\odot\bardxi
   +\barxi zd\xi\odot\bardz  +\eta\barz dz\odot\bardeta
   +\bareta z d\eta\odot\bardz 
   \right) \end{aligned}\eee
Each of the cross-terms in the last line can be estimated using the Cauchy
inequality (\ref{e3.4}), e.g.,
\bee
 \begin{aligned}
 |\xi\barz dz\odot\bardxi|
 &=|\xi dz\odot\overline{zd\xi}| \leq\frac{1}{2}
   (\xi dz\odot\overline{\xi dz} +zd\xi\odot\overline{zd\xi})\\
 &=\frac{1}{2} (|\xi|^2dz\odot\bardz+|z|^2d\xi\odot\bardxi).
 \end{aligned}\eee
Therefore
\bee
  \begin{aligned}
  \Phi^\ast g_{st}
  &\leq \frac{1}{2}((|\xi|^2+|\eta|^2) dz\odot\bardz 
    +(1+|z|^2)d\xi\odot\bardxi
    +(1+|z|^2)d\eta\odot\bardeta\\
  & +|\xi|^2dz\odot\bardz   +|z|^2d\xi\odot\bardxi
          +|\eta|^2dz\odot\bardz +|z|^2d\eta\odot\bardeta)\\
  =&\frac{1}{2}
    (2(|\xi|^2+|\eta|^2) dz\odot\bardz
   +(1+2|z|^2)d\xi\odot\bardxi
   +(1+2|z|^2)d\eta\odot\bardeta) \\
   &\hbox{$\leq 9g$ since $|z|<2$.}
  \end{aligned}\eee
To prove the inverse inequality let us go back to \refL{6}. There we
proved that
\bee
  \frac{1}{2} \frac{dz\odot\bardz}{(1+|z|^2)^2}
  \leq \frac{2}{|w|^2} \frac{1}{2}\Sum^4_{i=1} dw_i\odot\overline{dw_i},
  \eee
where $z=\frac{w_2}{w_1}=\frac{w_4}{w_3}$ when $w_1,w_3\not=0$.
Since $w=(\xi,z\xi,\eta,z\eta)$ in terms of $\xi,\eta$, this gives
\bee
  \begin{aligned}
 \frac{1}{2} \frac{dz\odot\bardz}{(1+|z|^2)^2}
    &\leq \frac{2}{|\xi|^2+|z\xi|^2+|\eta|^2+|z\eta|^2}\Phi^\ast g_{st},\\
 \text{and}\quad\Phi^\ast g_{st}
    &\geq\frac{1}{2}\ \frac{(|\xi|^2+|\eta|^2)}{2}
    \ \frac{dz\odot\bardz}{(1+|z|^2)}.
 \end{aligned}
\eee
Also obviously
\bee
  \Phi^\ast g_{st}\geq\frac{1}{2}(d\xi\odot\bardxi+d\eta\odot\bardeta)
\eee
Adding together the last two inequalities we obtain
\bee
  \begin{aligned}
  2\Phi^\ast g_{st}
   &\geq\frac{1}{2}(\frac{1}{2}(|\xi|^2+|\eta|^2)(1+|z|^2)^{-1}
    dz\odot\bardz+d\xi\odot\bardxi+d\eta\odot\bardeta) \\
  \Phi^\ast g_{st}
   &\geq\frac{1}{2\cdot 2\cdot 5}\cdot\frac{1}{2}
    (|\xi|^2+|\eta|^2)dz\odot\bardz+d\xi\odot\bardxi+d\eta\odot\bardeta)\\
   =&\frac{1}{20}g
  \end{aligned}\eee
Thus $20^{-1}g\leq\Phi^\ast g_{gs} \leq 20g$
and $\Phi$ is bi-Lipschitz.
\end{proof}

\begin{figure}\label{F:5}
\hskip2in\epsfig{file=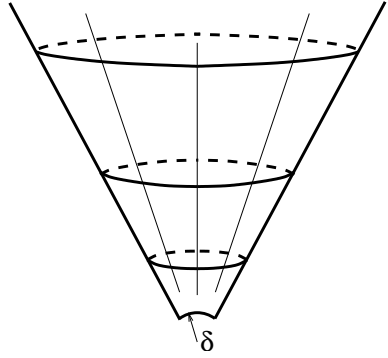,width=2truein}
\caption{Horizontal circles in the cone correspond to $\CP^1$-like sections of the conifold.}
\end{figure}

The geometric meaning of this lemma is best described by considering a
3-dimensional analog. Let $z^2=x^2+y^2$ be the standard cone in $\R^3$.
If $r,\theta$ are the polar coordinates in the $xy$ plane, then the induced
metric on the cone is $2dr^2+r^2d\theta^2$.
The metric $g$ on the conifold has similar structure with $\theta$ replaced
by $z$ and $r$ replaced by $\xi,\eta$. Sections with $r=\const$ describe
circles of growing radii as $r\to\infty$ (Fig.\,5). In a similar way sections
$|\xi|^2+|\eta|^2=\const$ describe disks of growing diameter along the fiber
parameters $\xi,\eta$.
Now recall that $2\D\times(\C\setminus\{0\})$ only parametrizes half
of the conifold, the other half is also $2\D\times(\C\setminus\{0\})$
with disks overlapping over an annulus containing the unit circle. Put together
the disks form $\CP^1\hbox{-like}$ sections of increasing area and diameter
as $\xi,\eta\to\infty$. 

On the narrow end, just as circles collapse
into the singular point of the usual cone at the origin, in the conifold
$\CP^1$ sections collapse into the conifold singularity when $\xi=\eta=0$.
In the resolved conifold this collapse is prevented by replacing the
singular point with a copy of $\CP^1$ and adding the Fubini-Study term
to the metric.

\begin{corollary}\label{C:4}
For any $\delta>0$ the second fundamental form of $\calC$ in $\C^4$ is uniformly bounded
on $\calC\setminus B_\delta(0)$.
\end{corollary}

\begin{proof}
If $w\in\calC\setminus\{0\}$ then $w=\Phi(z,\xi,\eta)$ by \refL{9}.
In fact, a whole neighborhood of $(z,\xi,\eta)$ is mapped into
a neighborhood of $w$ in $\calC$.
Therefore as  was commented after \refC{3}
$\|\II^{\calC/\C^4}(w)\|\leq C\|\nabla^2\Phi\|$,
where $\nabla^2\Phi$ is computed in $g_{st}$ on $\C^4$ and any metric
$g_1$ on $2\D\times(\C^2\setminus\{0\})$
such that $C^{-1}g_1\leq\Phi^\ast g_{st}$.
Let $g_1=\frac{1}{2}(dz\odot\bardz+d\xi\odot\bardxi+d\eta\odot\bardeta)$.
If $|w|^2=(1+|z|^2)(|\xi|^2+|\eta|^2)\geq \delta^2$, then
$(|\xi|^2+|\eta|^2)\geq\frac{\delta^2}{1+|z|^2}\geq\frac{\delta^2}{5}$
and by \refL{10}
\bee
  g_1\leq\frac{5}{\delta^2}g\leq\frac{5}{\delta^2}\cdot 20\,\Phi^\ast g_{st}
  =\frac{100}{\delta^2}\Phi^\ast g_{st}\eee
Thus one can take $C=\frac{100}{\delta^2}$. Both $g_1,g_{st}$ are flat.
Therefore by \refL{8} $\nabla^2\Phi=\partial^2\Phi$, the usual Hessian.
But the usual Hessian of $\Phi$ is in fact constant as one can see by
inspection (all entries are polynomials  in $z,\xi,\eta$ of at most
2nd degree). Therefore $\nabla^2\Phi$ and $\II^{\calC/\C^4}$ are
uniformly bounded on $\calC\setminus B_\delta(0)$.
\end{proof}

Since $C=\frac{100}{\delta^2}\underset{\delta\to\infty}{\longrightarrow} 0$
the second fundamental form is not
just bounded, it is asymptotically $0$. This means that the conifold
$\calC$ is asymptotically locally flat just like the usual cone.

We now want to extend this conclusion to the resolved conifold $\hatcalC$.
Recall that $\hatcalC\setminus 0(\hatcalC)$ can be obtained from
$\calC\setminus\{0\}$ by applying the biholomorphism $\pi^{-1}_2$.
To be able to use second fundamental forms in $\C^4$ and $\CP^1\times\C^4$
we need to extend $\pi^{-1}_2$ to $\C^4$.
Unfortunately this is impossible, but we can partially extend $\pi^{-1}_2$
into a neighborhood in $\C^4$ of every point $w\in\calC\setminus\{0\}$.
Namely, 
$$
\overline{\pi^{-1}_2}(w):=
\begin{cases} 
([w_1:w_2],w)&\text{if}\quad |w_1|^2+|w_2|^2\not=0\\                         
([w_3:w_4],w)&\text{if}\quad |w_3|^2+|w_4|^2\not=0.
\end{cases}
$$   
If $|w|^2>\delta^2$ then either $|w_1|^2+|w_2|^2>\frac{\delta^2}{2}$
or $|w_3|^2+|w_4|^2>\frac{\delta^2}{2}$ so one of the conditions is always satisfied and at least one of the extensions is defined around any point. Boundedness of $\II^{\hatcalC/(\CP^1\times\C^4)}$ follows from that of
$\II^{\calC/\C^4}$ by \refC{3} and the next Lemma. 

\begin{lemma}\label{L:11}
For any $\delta>0$ the extension $\overline{\pi^{-1}_2}$ is uniformly bi-Lipschitz
and has bounded covariant Hessian along $\calC\setminus B_\delta(0)$.
\end{lemma}

\begin{proof}
Due to symmetry it suffices to consider the case $|w_1|^2+|w_2|^2>\frac{\delta^2}{2}$ and
$\overline{\pi^{-1}_2}(w)=([w_1:w_2],w)\in\CP^1\times\C^4$.
Introduce the standard coordinate charts on $\CP^1$ with coordinate maps
$\varphi:[u:v]\mapsto\frac{u}{v}\in 2\D$ and $\mapsto\frac{v}{u}\in 2\D$.
Again due to symmetry it suffices to consider just one.
In coordinates we have
\bee
\begin{matrix}(\varphi,\id)\circ\overline{\pi^{-1}_2}: &\C^4\to 2\D\times\C^4\\
               &w\mapsto\left(\frac{w_1}{w_2},w\right)\end{matrix}.
\eee
Just as in \refL{6} one gets
\bee
  g_{st}\leq(\varphi\circ\overline{\pi^{-1}_2})^\ast \hatg
  \leq\left(1+\frac{2}{|w_1|^2+|w_2|^2}\right) g_{st}
  \leq\left(1+\frac{4}{\delta^2}\right)g_{st}.
\eee
Now we take a look at the covariant Hessian. By \refL{8}
\bee
  \nabla^2\Phi(X,Y)=\partial^2\Phi(X,Y)
  +\tildeG(\partial\Phi(X),\partial\Phi(Y))
  -\partial\Phi(\Gamma(X,Y)). 
\eee
In our case $\Phi=(\varphi,\id)\circ\overline{\pi^{-1}_2}$ and  $\Gamma=0$ since
$\C^4$ is flat. $\tildeG$ only depends on $z\in 2\D$ since
$\hatg=\pi^\ast_1 g_{FS}+\pi^\ast_2 g_{st}$ and
the second term is flat. Since $2\barD$ is compact and $\tildeG$
extends smoothly to $\C$, it is uniformly bounded on $2\D$.
Finally, since $\Phi$ is holomorphic we may consider
just holomorphic derivatives.

$\Phi$ is linear in $\C^4$ variables hence first derivatives are constant
and second ones are $0$. So the only part that matters is
$\Phi_1(w):=\varphi\circ\overline{\pi^{-1}_2}(w)=\frac{w_1}{w_2}$. By a direct computation:
$$
\partial\Phi_1=\left(\frac{1}{w_2},-\frac{w_1}{w^2_2}\right)
               =\left(\frac{1}{w_2},-\frac{z}{w_2}\right),
$$
where $z=\frac{w_1}{w_2}$ and
\bee
 \partial^2\Phi_1=
 \begin{pmatrix}0&-\frac{1}{w^2_2}\\
 -\frac{1}{w^2_2}&\frac{2w_1}{w^3_2}\end{pmatrix}
 =
 \begin{pmatrix}0&-\frac{1}{w^2_2}\\ 
 -\frac{1}{w^2_2}&\frac{2z}{w^2_2}\end{pmatrix}.
 \eee
Now recall that by our choice of coordinates $|z|<2$ so
$|w_2|>\frac{|w_1|}{2}$, $5|w_2|^2>|w_1|^2+|w_2|^2>\frac{\delta^2}{2}$ and hence 
$|w_2|>\frac{\delta}{\sqrt{10}}$. Therefore both $\partial\Phi_1$, $\partial^2\Phi_1$ 
and  therefore $\partial\Phi$, $\partial^2\Phi$ are uniformly bounded. Together with the previous
remarks this implies the same for $\nabla^2\Phi$.
\end{proof}
\begin{corollary}\label{C:5}
The resolved conifold $\hatcalC$ has bounded sectional curvature.
\end{corollary}

\begin{proof}
Note that $\pi^{-1}_2(\calC\bigcap\overline{B_\delta}(0))$ is a
compact subset in $\hatcalC$. Its complement in $\hatcalC$ is the
image under $\pi^{-1}_2$ of $\calC\setminus B_\delta(0)$.
By \refC{4} $\II^{\calC/\C^4}$ is uniformly bounded on
$\calC\setminus B_\delta(0)$. Every point in $\calC\setminus B_\delta(0)$
has a neighborhood in $\C^4$ such that $\pi^{-1}_2$ extends to it
and by \refL{11} these extensions are uniformly bi-Lipschitz with
bounded covariant Hessian. Therefore, by \refC{3},
$\II^{\hatcalC/(\CP^1\times\C^4)}$ is uniformly bounded on 
$\pi^{-1}_2(\calC\setminus B_\delta(0))$.
Since its complement has compact closure 
$\II^{\hatcalC/(\CP^1\times\C^4)}$ is bounded on the whole $\hatcalC$.
Since $\CP^1\times\C^4$ is a product of a closed manifold and a
flat manifold with the product metric $\sec^{\CP^1\times\C^4}$ is bounded.
Finally by the Gauss equation (\ref{e4.2}):
\bee
  \sec^\hatcalC(X,Y)=\sec^{\CP^1\times\C^4}(X,Y)
  +\hatg(\II^{\hatcalC}(X,X),\II^{\hatcalC}(Y,Y))
  -|\II^\hatcalC(X,Y)|^2_\hatg 
\eee
and all the terms on the right are bounded. Therefore so is $\sec^\hatcalC$.
\end{proof}

Now we want to establish that the injectivity radius $i(\hatcalC)$
is strictly positive. It is convenient to use the following criterion:
\begin{proposition}[Proposition 3.19 of \cite{ALN}]\label{P:3}
Let $(M,g)$ be a Riemannian manifold with complete metric and bounded
sectional curvature. Then three conditions are equivalent:

{\rm (i)} $i(M)>0$ ($i(M)$ is the injectivity radius.)

{\rm (ii)} There exist numbers $\delta$, $C>0$ such that every loop
$\gamma$ in $M$ of length $\l(\gamma)\leq\delta$ bounds a disc
$D$ in $M$ of $\diam(D)\leq C\cdot\l(\gamma)$.

{\rm (iii)} Every point in $M$ has a neighborhood uniformly bi-Lipschitz
to the flat unit ball.
\end{proposition}

In particular it follows directly from (iii) that

\begin{corollary}\label{C:6}
Let $\Phi:(M,g)\to(\tildeM,\tildeg)$ be a bi-Lipschitz diffeomorphism
between complete Riemannian manifolds with bounded sectional curvatures.
Then $i(M)>0$ if and only if $i(\tildeM)>0$.
\end{corollary}
$\hatcalC$ is obviously complete since it is properly embedded in
$\CP^1\times\C^4$ and it has bounded sectional curvature by
\refC{5}.

\begin{theorem}\label{T:3}
The resolved conifold $\hatcalC$ has bounded geometry.
\end{theorem}

\begin{proof}
Due to \refC{5} it only remains to prove that $i(\hatcalC)>0$.
From \refL{9} and the definition of $\hatcalC$ we have $\hatcalC$ covered
by the parametrizations:

\bee
  \begin{aligned}
  2\D\times\C^2
    &\overset{\Phi}{\lra}\hatcalC\hra\CP^1\times\C^4\\
  (z,\xi,\eta)
    &\longmapsto ([1:z],\xi,z\xi,\eta,z\eta);
  \end{aligned}\eee
\bee
  \begin{aligned}
  2\D\times\C^2
    &\overset{\Phi}{\lra}\hatcalC\hra\CP^1\times\C^4\\
   (z,\xi,\eta)
    &\longmapsto ([z:1],z\xi,\xi,z\eta,\eta);
  \end{aligned}\eee
By \refL{10} $\Phi^\ast\hatg$ is bi-Lipschitz to
\bee
  \begin{aligned}
  \Phi^\ast\pi^\ast_1 g_{FS}+g
  &=\frac{1}{2} \frac{dz\odot\bardz}{(1+|z|^2)} +\frac{1}{2}
   \left( (|\xi|^2+|\eta|^2)dz\odot\bardz+d\xi\odot\bardxi+d\eta\odot\bardeta\right) \\
  \text{and also to}\quad\tildeg &:=\frac{1}{2} \left((1+|\xi|^2+|\eta|^2)dz\odot\bardz
    +d\xi\odot\bardxi+d\eta\odot\bardeta \right).
  \end{aligned} \eee
Let $(z_0,\xi_0,\eta_0)\in 2\D\times\C^2$. We may assume
$|(\xi_0,\eta_0)|:=\sqrt{|\xi_0|^2+|\eta_0|^2}\geq 3$ since the
complement has compact closure. Consider the map
\bee
  \begin{aligned}
  \C^3\supset B_1(0) &\overset{f}{\lra}\C\times\C^2 \\
  (\lambda,\alpha,\beta) &\longmapsto \left(z_0+
   \frac{\lambda}{\sqrt{1+|\xi_0|^2+|\eta_0|^2}},\xi_0+\alpha,\eta_0+\beta\right).
  \end{aligned}
\eee
We claim that $f$ is uniformly bi-Lipschitz. Indeed,
\bee
  f^\ast\tildeg=\frac{1}{2} \left(
  \frac{1+|\xi|^2+|\eta|^2}{1+|\xi_0|^2+|\eta_0|^2}
   d\lambda\odot\bardlambda +d\alpha\odot\bardalpha +d\beta\odot\bardbeta
   \right)
\eee
To prove that $f^\ast\tildeg$ is equivalent to $g_{st}$ we need uniform estimates from
above and below. Note that by definition of $f$, $|(\xi,\eta)-(\xi_0,\eta_0)|\leq 1$.
Therefore,
\bee
  \begin{aligned}
  \frac{1+|\xi|^2+|\eta|^2}{1+|\xi_0|^2+|\eta_0|^2}=\frac{1+|(\xi,\eta)|^2}{1+|(\xi_0,\eta_0)|^2}
  & \leq
  \frac{1+(|\xi_0,\eta_0)|+|(\xi,\eta)-(\xi_0,\eta_0)|)^2}
     {1+|(\xi_0,\eta_0)|^2} \\
  & \leq
  \frac{(1+|(\xi_0,\eta_0)|^2) +2|(\xi_0,\eta_0)|+1}
     {1+|(\xi_0,\eta_0)|^2} \leq 4
  \end{aligned}\eee
On the other hand, 
\bee
  \begin{aligned}
  \frac{1+|\xi|^2+|\eta|^2}{1+|\xi_0,\eta_0|^2}
  & \geq
  \frac{(1+|(\xi_0,\eta_0)|^2) -2|(\xi_0,\eta_0)|)+1}
     {1+|(\xi_0,\eta_0)|^2}\\
  & =1-\frac{2|(\xi_0,\eta_0)|-1}{1+|(\xi_0,\eta_0)|^2}
    \geq \frac{1}{2} \hbox{\ when\ } |(\xi_0,\eta_0)|\geq 3.
  \end{aligned}\eee
Therefore regardless of the chosen point
\bee
  C^{-1}g_{st}\leq f^\ast\tildeg\leq C g_{st}\eee
with $C=4$. By \refP{3} we now have $i(\hatcalC)>0$ and $\hatcalC$
has bounded geometry.
\end{proof}

This proof is not very illuminating as to why $i(\hatcalC)>0$.
It is useful to have in mind the analogy between the conifold
and the usual cone in $\R^3$ described after \refL{10}. Let
$\calC_\delta=\calC\setminus B_\delta(0)$ and
$\hatcalC_\delta=\hatcalC\setminus\pi^{-1}_2(B_\delta(0))$.
Asymptotically geometries of $\calC_\delta$ and $\hatcalC_\delta$
are the same as one can see from expressions for $g$ and $\tildeg$
in the theorem because the Fubini-Study term becomes negligible as
$|\xi|^2+|\eta|^2\to\infty$. More precisely,
$\pi^{-1}_2:\calC_\delta\to\hatcalC_\delta$ is bi-Lipschitz with the
constant $C(\delta)\underset{\delta\to\infty}{\lra}1$.
Therefore `horizontal' sections of $\hatcalC$ are copies of $\CP^1$ with
K\"ahler volume $\underset{\delta\to\infty}{\lra}\infty$
just as horizontal sections of the cone are circles of increasing
diameters (see Fig.\,5). This means that cut points for points in $\hatcalC_\delta$
that are on the `other side' of the conifold are further and further
away from them as $\delta\to\infty$. Thus not only is $i(\hatcalC_\delta)$ bounded from below
but in fact $i(\hatcalC_\delta)\underset{\delta\to\infty}{\lra}\infty$.
Similarly as was noted after \refC{4},
$\sec(\calC_\delta)\underset{\delta\to\infty}{\lra}0$ and $\sec(\hatcalC_\delta)$ behaves 
the same way by the Gauss equation since $\pi^{-1}_2$ has bounded covariant Hessian. Summarizing, 
not only does $\hatcalC$ have bounded geometry but it is in fact asymptotically globally flat.

\section{Compactness of the moduli}\label{S:5}

In this section we prove the main result of this paper on moduli compactness of open pseudoholomorphic curves ending on the conifold transited perturbed conormal bundles. After briefly recalling the notions of open stable maps and their Gromov convergence we state the Sikorav compactness theorem (\refT{4}) and proceed to verify its assumtions in our case.

Let $(M,J,g)$ be an almost K\"ahler manifold,
$L\hra M$ be its totally real submanifold and $\Sigma$ be a Riemann surface with boundary $\partial\Sigma$ and a complex structure $j$. A smooth open pseudoholomorphic curve in $M$ ending on $L$ (or with the boundary on $L$) is a map $f:(\Sigma,\partial\Sigma)\to(M,L)$ such that $f_\ast j=Jf_\ast$.  
If instead of smooth Riemann surfaces one considers
complex one-dimensional varieties with at most nodal singularities
(i.e. stable curves with boundary) then $f$ is called an open
stable map \cite{L}. As is common in the literature, we often call stable maps pseudoholomorphic curves as well. Let $\partial\Sigma_j$ be the boundary
components of $\Sigma$, $\partial\Sigma=\cup_i\partial\Sigma_i$
and let $\alpha\in H_2(M,L)$ and $\beta_i\in H_1(L)$
be integral homology classes.

\begin{definition}\label{D:11}
The moduli space of open genus $g$ curves with $h$  boundary components is 
\begin{multline}
\overline{\calM}_{g,h}(M,L\mid \alpha,\beta_1,\dots,\beta_n):=\\
\{f\text{ open stable map}\mid\text{\,{\rm genus}\,}(\Sigma)=g,\#\{\partial\Sigma_i\}=h,f_\ast[\Sigma]=\alpha,
f_\ast[\partial\Sigma_i]=\beta_i\}. 
\end{multline}
\end{definition}

The appropriate topology on the moduli can be defined using the Gromov convergence. Since the domain of the limit curve may differ from that of the prelimit ones one needs some kind of smooth resolution of nodes to pull back the maps to the same domain. For open curves the definition of a resolution is worked out in \cite{L} where the interested reader is directed.

\begin{definition}[Gromov convergence]
One says that a sequence of stable maps $(\Sigma_n,f_n)$ Gromov converges to a map $(\Sigma,f)$ if there is a sequence of resolutions $\kappa_n:\Sigma_n\to\Sigma$ such that for any neighborhood $V$ of the union of all nodes in $\Sigma$:

{\rm 1)} $f_n\circ\kappa^{-1}_n\underset{n\to\infty}{\lra}f$ in $C^\infty(\Sigma\backslash V)$;

{\rm 2)} $(\kappa^{-1}_n)^\ast j_n\underset{n\to\infty}{\lra} j$ in $C^\infty(\Sigma\backslash V)$, where $j_n$, $j$ are complex structures on $\Sigma_n$, $\Sigma$ respectively;

{\rm 3)} $\Area(f_n(\Sigma_n))\underset{n\to\infty}{\lra}\Area(f(\Sigma))$.
\end{definition}

These moduli spaces have been used to define open Gromov-Witten
invariants in \cite{L, KL}. However, for such definitions to work
$\overline{\calM}_{g,h}$ at least has to be compact. Proving
compactness in the case $M=\hatcalC$ and $L=\CT(N^\ast_{k,\ve})$
will be our goal in this section. 

Compactness theorems for curves with boundary were considered by
several authors \cite{S,Ye}. For our purposes the most
suitable result is due to J.-C. Sikorav \cite{S}. We restate it here in a form consistent with the terminology used throughout the paper.

\begin{definition}\label{D:12}
$L\hra(M,J,g)$ satisfies the 2-point estimate if there exist constants
$C,\rho>0$ such that for any two points $x,y\in L$ with
$\dist^M(x,y)<\rho$ one has $\dist^L(x,y)\leq C\dist^M(x,y)$.
\end{definition}

\begin{theorem}[Proposition 5.1.2 and Theorem 5.2.3 of \cite{S}]\label{T:4}
Assume that $(M,J,g)$ has bounded geometry and $L\hra M$ is a
uniformly tame Lagrangian submanifold that satisfies the 2-point estimate.
Let $f_n:(\Sigma_n,\partial\Sigma_n)\to(M,L)$ be a sequence of open
curves with uniformly bounded areas such that
$f(\Sigma_n)\cap K\not=\emptyset$ for some compact subset $K\subset M$.
Then there exists a subsequence $f_{n_k}$ that Gromov-converges
to an open curve.
\end{theorem}
We will now verify the assumptions of the Sikorav theorem in the order they are listed.
Note that we already proved in \refS{4} that $\hatcalC$ is geometrically
bounded and in \refS{3} that $L=\CT(N^\ast_{k,\ve})$ is a tame Lagrangian in it.

The 2-point estimate is preserved under bi-Lipschitz maps,
i.e., if $L\hra M$ satisfies it and $f:M\to\tildeM$ is
bi-Lipschitz then so does $f(L)\hra\tildeM$.
Moreover, locally any submanifold satisfies it.

\begin{lemma}\label{L:12}
Let $L\hra(M,g)$. Then every point of $L$ has a neighborhood
in $L$ such that $\dist^L(q,q')\leq C\dist^M(q,q')$ for
$q,q'$ in this neighborhood for some constant $C$.
\end{lemma}

\begin{proof}
Let $q_0\in L$. Choose a neighborhood $U_{q_0}$ of $q_0$
in $M$ with coordinate function $\varphi$ such that $\varphi(q_0)=0$,
$\varphi(L\cap U_{q_0})=\{x\in\R^m\mid x_{\l+1}=\dots=x_m=0\}$,
where $m:=\dim M$, $\l:=\dim L$.
For any ball $B_R\subset\R^m$ the map $\varphi:\varphi^{-1}(B_R)\to(\R^m,g_{st})$
is bi-Lipschitz by a compactness argument. Fix $R$ and let $C$ be the
corresponding bi-Lipschitz constant. Then
\bee
  \begin{aligned}
  \dist^L(q,q')
  &\leq C^{1/2}\dist_{\R^l}(\varphi(q),\varphi(q'))
     =C^{1/2}\dist_{\R^m}(\varphi(q),\varphi(q'))\\
  &\leq C^{1/2}\cdot C^{1/2}\dist^M(q,q') = C\dist^M(q,q')
  \end{aligned}\eee
\end{proof}

\begin{figure}\label{F:6}
\hskip2.5in\epsfig{file=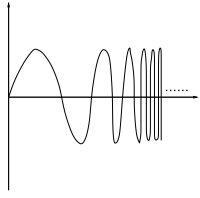,width=1truein}
\caption{Submanifold in $\R^2$ that does not satisfy the 2-point estimate.}
\end{figure}

This implies of course that any compact submanifold satisfies the
2-point estimate. An example of a submanifold that does not satisfy it
is given by the graph of $\sin(\pi e^x)$ in $\R^2$. This graph is the the graph of sine that gets more and more compressed as $x\to\infty$ (see Fig.\,6). The distance in $\R^2$ between two consecutive zeros is $\ln(n+1)-\ln(n)\to0$. However the distance between them along the graph is $\geq2$ since one has to go along the arc of sine to get from one to the next.

Recall that $\CT(N^\ast_{k,\ve})=\pi^{-1}_2\circ\Phi_\ve(F_\ve(N^\ast_k))$,
where $\Phi_\ve(x,p)=(x,p+\ve\xi(x))$ and
$F_\ve(x,p)=(x\sqrt{|p|^2+\ve^2 },p)$.
Since both $\Phi_\ve$ and $\pi^{-1}_2$ are uniformly bi-Lipschitz
on the relevant sets it suffices to prove the 2-point estimate for
$F_\ve(N^\ast_k)=\{(k(t)\sqrt{|p|^2+\ve^2},p)\in\R^4\times\R^4
\mid p\cdot k(t)=p\cdot\dotk(t)=0,t\in S^1\}$.

\begin{lemma}\label{L:13}
$F_\ve(N^\ast_k)\hra\R^4\times\R^4$ and therefore
$\CT(N^\ast_{k,\ve})\hra\CP^1\times\C^4$ satisfy the 2-point estimate.
\end{lemma}

\begin{proof}
Despite the length this proof reduces to multiple applications of the triangle inequality. For convenience, it is split into three steps corresponding to pairs of points at a fixed distance from the zero section, at a distance greater than $1$ from it and finally in the general position.

{\sc\underline{ Step 1 }} Let $S_r:=\{(x,p)\in F_\ve(N^\ast_k)\mid |p|=r\}$. In this step we will prove that they 
satisfy the 2-point estimate in $\R^4\times\R^4$ with constants
$C,\rho$ independent of $r\geq 1$.

Since $S_1$ is compact in $\R^4\times\R^4$, by \refL{12} it satisfies
the estimate with some $C_1,\rho_1>0$. We claim that the same
$\rho_1$ and $C_1\sqrt{1+\ve^2}$ work for all $r\geq 1$.
Indeed, let $q=(x,p)$, $q'=(x',p')\in S_r$. Then
$q_1:=\left(\frac{x}{\sqrt{r^2+\ve^2}},\frac{p}{r}\right)$,
$q'_1:=\left(\frac{x'}{\sqrt{r^2+\ve^2}},\frac{p'}{r}\right)\in S_1$.
Since
\bee
  \dist^M(q_1,q'_1)
  =\sqrt{\frac{1}{r^2+\ve^2}|x-x'|^2+\frac{1}{r^2}|p-p'|^2} 
  \leq \dist^M(q_1,q'_1)<\rho_1, \eee
where $M=\R^4\times\R^4$, we have
$\dist^{S_1}(q_1,q)\leq C\cdot\dist^M(q_1,q'_1)$.

Let $\gamma=(\gamma_x,\gamma_p)$ be any path in $S_1$ connecting
$q_1$ and $q'_1$. Then $\gamma^r:=(\gamma_x\sqrt{r^2+\ve^2},r\gamma_p)$ is a path in $S_r$ connecting $q$ and $q'$. Moreover, $r\l(\gamma)\leq\l(\gamma^r)\leq\sqrt{r^2+\ve^2}\,\l(\gamma)$.
Minimizing over all such paths one gets
\bee
  r\dist^{S_1}(q_1,q'_1)\leq\dist^{S_r}(q,q')
  \leq\sqrt{r^2+\ve^2}\dist^{S_1}(q_1,q'_1) \eee
Also obviously $r\dist^M(q_1,q'_1)\leq\dist^M(q,q')\leq\sqrt{r^2+\ve^2}\dist^M(q_1,q'_1)$.
Thus
\bee
  \begin{aligned}
  \dist^{S_r}(q,q')
  &\leq\sqrt{r^2+\ve^2}\dist^{S_1}(q_1,q'_1)
   \leq C\sqrt{r^2+\ve^2}\dist^M(q_1,q'_1)\\
  &\leq\frac{C\sqrt{r^2+\ve^2}}{r}\dist^M(q,q')
   \leq C_1\sqrt{1+\ve^2}\dist^M(q,q') 
  \end{aligned}\eee
and $C_1\sqrt{1+\ve^2},\rho_1$ work for all $r\geq 1$.

\begin{figure}\label{F:7}
\hskip1.8in\epsfig{file=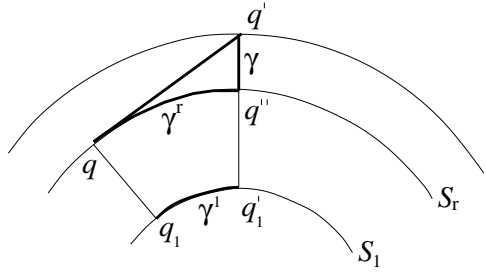,width=2.5truein}
\caption{The 2-point estimate.}
\end{figure}

\vspace{.15in}
\nd{\sc\underline{ Step 2 }}
Let $q:=(x\sqrt{|p|^2+\ve^2},p)$,
$q':=(x'\sqrt{|p'|^2+\ve^2},p')\in F_\ve(N^\ast_k)$
and $1\leq|p|\leq|p'|$. Then the 2-point estimate holds with
$\tildeC_1=\sqrt{2}+3C_1\sqrt{1+\ve^2}$, $\tildep_1=\frac{\rho_1}{3}$. A proof follows.

Define $q'':=\left(x'\sqrt{|p|^2+\ve^2},|p|\frac{p'}{|p'|}\right)$,
then $q,q''\in S_{|p|}$ (see Fig.\,7).
Consider the following path:
\bee
  \gamma(t)=(x'\sqrt{t^2|p'|^2+\ve^2},tp'), \eee
then $\gamma(1)=q'$, $\gamma\left(\frac{|p|}{|p'|}\right)=q''$
and $\gamma\left([\frac{|p|}{|p'|},1]\right) \subset F_\ve(N^\ast_k)$.
Essentially, $q',q''$ are on the same radial line and $\gamma$
is the segment of the line connecting them.
\bee
  \begin{aligned}
  \dotg(t)
    &=\left(\frac{tx'}{\sqrt{t^2|p'|^2+\ve^2}},p'\right) \\
  \l(\gamma)
    &=\Int^1_{\frac{|p|}{|p'|}} \sqrt{g_{st}(\dotg,\dotg)} dt
     =\Int^1_{\frac{|p|}{|p'|}} \sqrt{\frac{t^2}{t^2|p'|^2+\ve^2}+|p'|^2}\,dt\\
&\leq \Int^1_{\frac{|p|}{|p'|}} \sqrt{1+|p'|^2}\,dt
     \leq\sqrt{2} |p'| |1-\frac{|p|}{|p|'}|
     \leq\sqrt{2} |p-p'|
     \leq\sqrt{2}\dist^M(q,q').
  \end{aligned}\eee
Now we want to estimate $\dist^M(q,q'')$ in terms of $\dist^M(q,q')$.
\bee
  \dist^M(q,q'')=
  \sqrt{(|p|^2+\ve^2)|x-x'|^2+|p-|p|\frac{p'}{|p'|}|^2}. \eee
Let us consider separately each term under the square root:
\bee
  \begin{aligned}
  |x-x'|\sqrt{|p|^2+\ve^2} 
  &\leq |x\sqrt{|p|^2+\ve^2}-x'\sqrt{|p'|^2+\ve^2}|
     + |x\sqrt{|p'|^2+\ve^2}-x'\sqrt{|p|^2+\ve^2}|\\
  &\leq\dist^M(q,q')
     + \frac{|(|p| +|p'|)(|p|-|p'|)}{\sqrt{|p'|^2+\ve^2}+\sqrt{|p|^2+\ve}}\\
  &\leq\dist^M(q,q') +|p-p'|\leq2\dist^M(q,q');\\
  |p-|p|\frac{p'}{|p'|}| 
  & \leq\frac{| p|p'|-|p|p'|}{|p'|}
     +\frac{| |p'|(p-p')+p'(|p'|-|p|)|}{|p'|}\\
  &\leq\frac{2|p'| |p-p'|}{|p'|} \leq 2\dist^M(q,q').
  \end{aligned}\eee
Therefore: $\dist^M(q,q'')\leq\sqrt{4+4}\dist^M(q,q')\leq 3\dist^M(q,q')$.

Finally, let $\dist^M(q,q')<\frac{\rho_1}{3}$, then $\dist^M(q,q'')<\rho_1$
and by the triangle inequality and Step 1 $(L=F_\ve(N^\ast_k))$:
\bee
  \begin{aligned}
  \dist^L(q,q')
  & \leq \dist^L(q,q'')+\dist^L(q'',q')\\
  & \leq C_1\sqrt{1+\ve^2} \dist^M(q,q'')+\l(\gamma)\\
  & \leq 3C_1\sqrt{1+\ve^2} \dist^M(q,q')+\sqrt{2}\dist^M(q,q')\\
  &=(3C_1\sqrt{1+\ve^2}+\sqrt{2}) \dist^M(q,q').
  \end{aligned}\eee

\vspace{.15in}
\nd{\sc\underline{ Step 3 }}
Now let $q,q'\in F_\ve(N^\ast_k)$ be two arbitrary points with
$\dist^M(q,q')<\frac{\rho_1}{3}$.
If $|p|,|p'|\geq 1$ they are covered by Step 2. Otherwise, both
are contained in
$\{(x,p)\in F_\ve(N^\ast_k)\mid |p|\leq 1+\frac{\rho_1}{3}\}$.
This set is compact and can be covered by a finite number of neighborhoods
as in \refL{12}. So there are $C_2,\rho_2$ that realize the 2-point
estimate there. Set $\rho:=\min\left(\rho_2,\frac{\rho_1}{3}\right)$,
$C:=\max(C_2,3C_1\sqrt{1+\ve^2}+\sqrt{2})$.  The entire $F_\ve(N^\ast_k)$
now satisfies the 2-point estimate with these $C,\rho$.
\end{proof}

\begin{corollary}[2-point estimate]\label{C:7} 
$\CT(N^\ast_{k,\ve})$ in $\hatcalC$ satisfies the 2-point estimate.
\end{corollary}

\begin{proof}
By \refL{13} this is true for $\CT(N^\ast_{k,\ve})$ in $\CP^1\times\C^4$.
Let $C,\rho$ be the constants. If
$\dist^{\hatcalC}(q,q')<\rho$ then $\dist^{\CP^1\times\C^4}(q,q')<\rho$
and $\dist^{\CT(N^\ast_{k,\ve})}(q,q')\leq C\dist^{\CP^1\times\C^4}(q,q')
\leq C\dist^{\hatcalC}(q,q')$ so the same constants work.
\end{proof}

The next step in meeting the assumptions of \refT{4} is to establish an
area bound for curves in the moduli space. When the submanifold they are
ending on is Lagrangian with respect to a K\"ahler form on the ambient manifold any two curves
in the same relative homology class have the same area.
Indeed, let $S$ be a chain realizing the relative homology.
Then $\partial S=\Sigma_1-\Sigma_2+\partial S\cap L$ and
$$
0=\int_Sd\omega=\int_{\partial S}\omega=\int_{\Sigma_1}\omega-\int_{\Sigma_2}\omega
+\int_{\partial S\cap L}\omega=\int_{\Sigma_1}\omega-\int_{\Sigma_2}\omega
$$
as $\omega|_L=0$.
But for pseudoholomorphic curves $\Area(\Sigma)=\int_\Sigma\omega$
so $\Area(\Sigma_1)=\Area(\Sigma_2)$. In our case we only have a symplectic form
$\tildeo:=\tildeo_\ve$ defined on $\hatcalC\setminus 0(\hatcalC)$ and uniformly tame
on every
$\hatcalC_\delta:=\hatcalC\setminus\pi^{-1}_2(\calC\cap B_\delta(0))$
that vanishes on $\CT(N^\ast_{k,\ve})$ (\refT{2}).

\begin{lemma}[Area Bound]\label{L:14} 
Let $L=\CT(N^\ast_{k,\ve})$, $\beta\in H_2(\hatcalC,L)$.
There exists a constant $A_\beta$ such that if $\Sigma$ is an open curve ending
on $L$ with $[\Sigma]=\beta$ then $\Area(\Sigma)\leq A_\beta$.
\end{lemma}

\begin{proof}
Define $\tildeg(X,Y):=\frac{1}{2}(\tildeo(X,JY)+\tildeo(Y,JX))$.
Then it follows from the tameness condition for $\tildeo$ that $\tildeg$
is a metric on $\hatcalC\setminus 0(\hatcalC)$ equivalent to $\hatg$
on every $\hatcalC_\delta$. Let $\tildeA$ denote the area with respect
to this metric. Then $C(\delta)^{-1}\Area\leq\tildeA\leq C(\delta)\Area$
by equivalence of metrics for surfaces in $\hatcalC_\delta$.
Just as in the case of compatible forms one proves that if
$\Sigma$ is pseudoholomorphic then $\tildeA(\Sigma)=\int_\Sigma\tildeo$
(see, e.g., \cite{MS}).

Recall from the discussion after \refC{1} that
$L=\CT(N^\ast_{k,\ve})\subset\hatcalC_{2\ve}$.
Let $\Sigma_1,\Sigma_2$ be two open curves ending on $L$ with
$[\Sigma_1]=[\Sigma_2]=\beta$ and let $S$ realize the relative homology,
i.e., $[\partial S]=[\Sigma_1]-[\Sigma_2]\mod L$.
Since $S$ is 3-dimensional, $0(\hatcalC)\simeq\CP^1$ is
2-dimensional and $\hatcalC$ is 6-dimensional 
we can push $S$ out of an $\ve\hbox{-neighborhood}$ of $0(\hatcalC)$ by transversality.
In other words, without loss of generality $S\subset\hatcalC_\ve$.
Since $\tildeo|_L=0$ just as above:
\bee
  0=\int_Sd\tildeo=\int_{\partial S}\tildeo
  =\int_{\Sigma_1}\tildeo-\int_{\Sigma_2}\tildeo
  =\tildeA(\Sigma_1)-\tildeA(\Sigma_2), \eee
and $\tildeA(\beta):=\tildeA(\Sigma)$, $[\Sigma]=\beta$ only depends on the
homology class. But by equivalence of metrics,
$\Area(\Sigma)\leq C(\ve)\tildeA(\beta)=: A_\beta$ for any
$[\Sigma]=\beta$.
\end{proof}

For pseudoholomorphic curves an area bound also implies a diameter
bound due to the monotonicity lemma. In particular, we have the following
property:

\begin{lemma}[Area-diameter estimate, Proposition 4.4.1 of \cite{S}]
\label{L:15}
Let $(M,J,g)$ have bounded geometry, $L\hra M$ be uniformly tame
Lagrangian with taming constant $C_1$ and satisfy the 2-point estimate with
constants $C_2,\rho$. Let $K\subset M$ be a compact subset and $\Sigma$
a pseudoholomorphic curve, $\partial\Sigma\subset L$,
$\Sigma\cap K\not=\emptyset$. Then $\Sigma\subset B_R(K)$
with $R=C\Area(\Sigma)$, $C:=\frac{8(C_1+C_2+1)}{\pi\min(i(M),\rho)}$.
In particular, $\diam(\Sigma)\leq C\cdot \Area(\Sigma)$.
\end{lemma}

\begin{proof}
Actually, in \cite{S} this is proved for closed curves, but the same proof
works for open ones if one uses the monotonicity lemma of Proposition 4.7.2
instead of that of Proposition 4.3.1. As for the second statement, it
suffices to take $K=\{\sigma\}$, where $\sigma\in\Sigma$ is any point.
\end{proof}

Geometrically this estimate means that pseudoholomorphic
curves can not be long and thin. But even when all assumptions of the last lemma are satisfied it does not
follow that any sequence of curves with bounded area has a convergent
subsequence.

\begin{figure}\label{F:8}
\hskip2.5in\epsfig{file=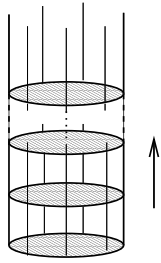,width=.8truein}
\caption{A sequence of holomorphic disks without Gromov-convergent subsequences.}
\end{figure}

\begin{example}\label{3}
Let $M=\C^2$, $L=S^1\times\R\subset\C\times\C=\C^2$,
$S^1:=\{z\in\C\mid |z|=1\}$.
One can easily check that $L$ is a Lagrangian submanifold with respect
to the standard symplectic form on $\C^2$ that satisfies the 2-point estimate
with $\rho=\infty$ and $C=\frac{\pi}{2}$.
$M$ of course has bounded geometry with $\sec^M=0$, $i(M)=\infty$.
Let $\Sigma=\barD=\{z\in\C\mid |z|\leq 1\}$ and
$f_n:\Sigma\to M$, $f_n(z)=(z,n)$.
$\Area(f_n)=\pi$ for all $n$, $\partial f_n(\Sigma)\subset L$, but
$f_n$ clearly does not have Gromov convergent subsequences (Fig.\,8).

The additional property which is missing here is some kind of
convexity of $L$ at $\infty$ that would force pseudoholomorphic
curves to be anchored to a compact subset.
In terms of \refT{4}, there is no compact subset in $\C^2$ that
all $f_n(\Sigma)$ meet.
\end{example}
In our case, anchoring to a compact subset follows from a simple observation below that generalizes the vanishing theorem of Witten (see \refR{1}). Despite the simplicity it implies that all relevant pseudoholomorphic curves in the resolved conifold must intersect its zero section.
\begin{lemma}[Anchoring]\label{L:16}
Let $M$ be a manifold with an almost complex structure $J$ and an exact $2$-form $\omega=d\lambda$ which is symplectic and tames $J$ on $M\backslash Z$. If $L\subset M\backslash Z$ is an exact Lagrangian submanifold then any non-constant pseudoholomorphic curve, closed or ending on $L$, intersects $Z$.
\end{lemma}
\begin{proof}
Suppose not, then $f(\Sigma)\subset M\backslash Z$. Let $g$ be the metric 
on $M\backslash Z$ determined by $\omega$ and $J$ as in \refL{14} and consider the corresponding area. Then as in \refR{1}
\be
\Area(f)=\int_\Sigma f^\ast\omega=\int_\Sigma d(f^\ast \lambda) =\int_{\partial\Sigma} f^\ast\lambda=[\lambda\!\mid_L](f_\ast[\partial\Sigma])=0. 
\ee
and we arrive at a contradiction with non-constancy of $f$.
\end{proof}
Note that if $Z=\emptyset$, i.e. $M$ is exact symplectic the lemma implies that it has no non-constant pseudoholomorphic curves as observed by Witten for $M$ a cotangent bundle \cite{W92}. In our case $M=\hatcalC$, $Z=0(\hatcalC)$ is the zero section, $L=\CT(N^\ast_{k,\ve})$ and $\omega=\tildeo_\ve$ is supplied by \refT{2}. Note that here $Z\simeq\CP^1$ and $f(\Sigma)$   are 2-dimensional and $\hatcalC$ is 6-dimensional. Thus the non-empty intersection granted by \refL{16} is a purely symplectic phenomenon that does not follow from a dimension count in differential topology. Now we are ready to prove the main result.

\begin{theorem}\label{T:5}
Moduli spaces
$\overline{\calM}_{g,h}(\hatcalC,\CT(N^\ast_{k,\ve})
 \mid\alpha,\beta_1,\dots,\beta_n)$
are compact for $\alpha\not=0\in H_2(\hatcalC,\CT(N^\ast_{k,\ve}))$ ($\alpha=0$ corresponds to constant maps).
\end{theorem}

\begin{proof}
Let $f_n$ be a sequence of open curves in the above moduli space. By \refL{13} $\CT(N^\ast_{k,\ve})$ satisfies the 2-point estimate, by \refL{14} $\Area(f_n)\leq A_\alpha<\infty$ and by \refL{16} all $f_n$ intersect $0(\hatcalC)\simeq\CP^1$ which is compact. Therefore by \refT{4} there exists a Gromov-convergent subsequence.  
\end{proof}
As a bonus, \refL{15} yields that for any pseudoholomorphic curve in the above moduli space $\diam(f)\leq C\cdot A_\alpha$ 
and we have them all contained in the ball $B_{C\cdot A_\alpha}(0(\hatcalC))$.

\section*{Conclusions}\label{S:6}

Whether the submanifolds constructed in this paper are suitable for the Gopakumar-Vafa conjecture remains to be seen. 
In addition to compactness one also needs the open curves moduli to have virtual dimension zero \cite{KL}. For open curves ending on a Lagrangian submanifold this is traditionally ensured by the special Lagrangian condition \cite{J}: the imaginary part of the holomorphic volume form vanishes along the submanifold. Since the holomorphic volume form on a Calabi-Yau is nowhere zero this means in particular, that its real part does not vanish along the submanifold. The special Lagrangian condition is far from being necessary. It suffices to have it satisfied only cohomologically: the Maslov class \cite{Ar} of the submanifold should be trivial. 

The last condition makes sense for totally real submanifolds as well as for Lagrangian ones. It holds for example if the real part of the holomorphic volume form does not vanish. Away from the zero section the holomorphic volume form for the resolved conifold can be obtained by deforming the one for the cotangent bundle. One can see that the real part of the latter is uniformly separated from zero along the conormal bundles to knots. For this reason we beleive that non-vanishing along their conifold transitions can be proved by a perturbation argument akin to the one used for the symplectic form in \refS{3}. 

Assuming that the Maslov class is zero one still needs to compare the Gromov-Witten invariants to the Chern-Simons knot invariants. A computational approach used so successfully for comparing the closed invariants \cite{FP,GR} does not seem to be feasible beyond the case of the unknot \cite{KL}. The problem is that the circle symmetry of the unknot is lacking in general and the standard localization techniques do not apply. It seems more likely that a proof will come from a deformation argument relating the Gromov-Witten theory on the resolved conifold to a degenerate string theory on the cotangent bundle postulated by Witten \cite{W92}. The recent work on contact knot homology \cite{Ng} is an interesting step in this direction.

{

}

\end{document}